\documentclass[a4paper,10pt]{amsart}
\usepackage[utf8]{inputenc}
\usepackage{lmodern}
\usepackage[T1]{fontenc}
\usepackage{microtype} 	
\usepackage[dvipsnames]{xcolor}
\usepackage[all]{xy}
\usepackage{amsmath,amssymb,amsfonts,amsthm,latexsym,mathrsfs}
\usepackage{ stmaryrd } 
\usepackage{tikz-cd}
\usepackage{hyperref}
\hypersetup{colorlinks=true,linkcolor=Maroon,citecolor=NavyBlue}
\theoremstyle{plain}
\newtheorem{theorem}[subsubsection]{{\bf Theorem}}
\newtheorem*{theorem*}{{\bf Theorem}}
\newtheorem{corollary}[subsubsection]{{\bf Corollary}}
\newtheorem*{corollary*}{{\bf Corollary}}
\newtheorem{proposition}[subsubsection]{{\bf Proposition}}
\newtheorem{lemma}[subsubsection]{{\bf Lemma}}

\theoremstyle{definition}

\theoremstyle{remark}
\newtheorem{remark}[subsubsection]{{\it Remark}}
\newtheorem{example}[subsubsection]{{\it Example}}
\numberwithin{equation}{section}
\DeclareMathOperator{\B}{B}
\DeclareMathOperator{\Br}{Br}
\DeclareMathOperator{\Tors}{Tors}
\DeclareMathOperator{\K}{K}
\DeclareMathOperator{\SK}{SK}
\DeclareMathOperator{\dd}{d}
\DeclareMathOperator{\HH}{H}
\DeclareMathOperator{\im}{im}
\DeclareMathOperator{\codom}{codom}
\DeclareMathOperator{\dom}{dom}

\DeclareMathOperator{\M}{M}
\DeclareMathOperator{\GL}{GL}
\DeclareMathOperator{\Gal}{Gal}

\DeclareMathOperator{\pr}{pr}

\DeclareMathOperator{\Spec}{Spec}
\DeclareMathOperator{\Proj}{Proj}
\DeclareMathOperator{\Sym}{Sym}
\newcommand{\FF}{\mathbb{F}}
\newcommand{\PP}{\mathbb{P}}
\newcommand{\Gr}{\mathbf{Gr}}

\newcommand{\QQ}{\mathbb{Q}}
\newcommand{\ZZ}{\mathbb{Z}}
\newcommand{\CC}{\mathbb{C}}
\newcommand{\Lie}{\mathcal{L}}
\newcommand{\II}{\mathcal{I}}

\newcommand{\pwedge}{\hat{\wedge}}
\newcommand{\grasslie}{{\bf \Lambda}}
\newcommand{\der}{\mathrm{d}}

\usepackage{enumitem}
\begin{document}
\baselineskip=15pt
\title[Irrationality of generic quotient varieties]{Irrationality
of generic quotient varieties \\ via Bogomolov multipliers}

\author[U. Jezernik]{Urban Jezernik}
\address{Faculty of Mathematics and Physics, University of Ljubljana and Institute of Mathematics, Physics, and Mechanics, Jadranska 19, 1000 Ljubljana, Slovenia}
\email{urban.jezernik@fmf.uni-lj.si}

\author[J. Sánchez]{Jonatan Sánchez}
\address{Department of Applied Mathematics (DMATIC), ETSI Ingenieros Informáticos, Uni- versidad Politécnica de Madrid, Campus de Montegancedo, Avenida de Montepríncipe, 28660, Boadilla del Monte, Spain}
\email{jonatan.sanchez@upm.es}

\thanks{
	The first author has received funding from the European Union’s Horizon 2020
research and innovation programme under the Marie Sklodowska-Curie grant
agreement No 748129. He has also been supported by the Spanish Government
grant MTM2017-86802-P and by the Basque Government grant IT974-16. The second
author was supported by the Spanish Government grant MTM2014-54804-P and the
Basque Government grant IT1094-16.
}
\date{\today}
\begin{abstract}
\noindent
The Bogomolov multiplier of a group is the unramified Brauer group
associated to the quotient variety of a faithful representation of the group.
This object is an obstruction for the quotient variety to be stably rational.
The purpose of this paper is to study these multipliers associated to
nilpotent pro-$p$ groups by
transporting them to their associated Lie algebras. Special focus is set on
the case of $p$-adic Lie groups of nilpotency class $2$, 
where we analyse the moduli space.
This is then applied to give information on asymptotic behaviour of 
multipliers of finite images of such groups of exponent $p$. We show that with fixed $n$ and
increasing $p$, a positive proportion of these groups of order $p^n$
have trivial multipliers. On the other hand, we show that by fixing $p$
and increasing $n$, log-generic groups of order $p^n$ have non-trivial
multipliers. Whence quotient varieties of faithful representations of 
log-generic $p$-groups are not stably rational.
Applications in non-commutative Iwasawa theory are developed.

\end{abstract}
\maketitle
\section{Introduction}

\subsection{The rationality problem}

Let $X$ be a smooth connected projective complex variety.
The famous {\em rationality problem} asks whether or not $X$ is 
birational to a projective space. 
In terms of function fields, this means that its
field of rational functions $\CC(X)$ is purely
transcendental over $\CC$.
This problem is especially interesting if one assumes that $X$ is 
a priori {\em unirational}, meaning that there is an inclusion
of $\CC(X)$ into some $\CC(t_1, \dots, t_n)$.
Under this additional assumption, the rationality problem
is known as the {\em Lüroth problem}. In its stable form,
it asks whether unirational varieties are {\em stably rational},
meaning that some purely transcendental extension of $\CC(X)$
is purely transcendental over $\CC$.

\subsection{Counterexamples to the Lüroth problem}

The answer to the Lüroth problem turns out to be positive for
curves and surfaces, and negative in general.
The first counterexamples were constructed 
independently by Clemens and Griffiths \cite{CleGri72}, Iskovskikh and Manin
\cite{IskMan71},
and by Artin and Mumford \cite{AM72}. The latter even constructed counterexamples
to the stable version of the problem. This is done by showing
that the {\em Brauer group},
which in this situation coincides with
$\Br(X) = \Tors \HH^3_{sing}(X_{an}, \ZZ)$,
of the specific $3$-folds they consider is non-trivial, and this is an
obstruction for such a variety to be stably rational.
Many more counterexamples have recently been produced
relying on the much finer intermediate property of possessing a  
{\em decomposition of the diagonal of the Chow group},
providing means to show that very general varieties in many
interesting families are not stably rational
(see \cite{Voi15}). 
The latter ultimately leads to a far-reaching generalization of the
counterexamples of Artin and Mumford.

\subsection{Unramified Brauer group}

In this paper, we address questions of a similar flavour as
the very general celebrated result mentioned above, but we stay with
the Brauer group as an obstruction to stable rationality. 
This obstruction can be made to work even for varieties
that are not smooth by passing to a smooth birational model
and computing the Brauer group there. The reason for this
is that the Brauer group of a smooth variety can be expressed as
the subgroup of the classical $\Br(\CC(X)) = \HH^2(\Gal(\CC(X)), \mathbb{G}_m)$
consisting of the cocycles
\[
\Br(X) = 
\bigcap_{\text{$A$ DVR of $\CC(X)$}} 
\ker 
\left( 
\Br(\CC(X))
\xrightarrow{\partial_A}
\HH^1(\kappa_A, \QQ/\ZZ)
\right),
\]
where $\partial_A$ is the residue map
(see \cite{GilSza17}).
This description of the Brauer group depends only on the function field
$\CC(X)$ and so it can be computed in any birational model of $X$.
The group on the right hand side is the {\em unramified Brauer group}
of $X$, denoted by $\Br_{nr}(X)$.

\subsection{Quotients by group actions}

Throughout the paper, we will deal with varieties $X$ that are GIT quotients
of affine spaces $V$ by faithful linear actions of finite groups $G$
(see \cite{ColSan07}).
Thus $\CC(X) = \CC(V)^G$. 
The Lüroth problem for such varieties is known as {\em Noether's problem},
especially when considered over $\QQ$ rather than $\CC$,
and it is tightly related to a constructive approach to the inverse
Galois problem (see \cite{JLY02}).
In the setting of quotients by group actions, Bogomolov was able to transfer
the above formula for the unramified Brauer group into purely group-theoretical
cohomology (see \cite{Bog87}),
\[
\Br_{nr}(\CC(V)^G)
\cong
\bigcap_{A \leq G, \, [A,A] = 1} \ker \left( 
\HH^2(G, \QQ/\ZZ)
\xrightarrow{\text{$res^G_A$}}
\HH^2(A, \QQ/\ZZ)
\right).
\]
The group on the right hand side is hence known in the literature as the
{\em Bogomolov multiplier}, denoted by $\B_0(G)$.
This description of the unramified Brauer group as the Bogomolov
multiplier has made it much more tractable to compute these invariants
and thus produce new counterexamples to rationality
as well as understand these objects better (see \cite{Mor12,HKK13,JezMor15unicp}).
It also explains that, in a way, the simplest obstruction for a
quotient variety $V/G$ to be stably rational is that the group $G$ possesses 
some special $2$-cocycles, a property that might be seen as rather special.

\subsection{Profinite groups}

The notion of the Bogomolov multiplier can be extended in a natural way
to the category of all profinite groups using the same cohomological definition.
Much work has been done in understanding Bogomolov multipliers of absolute Galois
groups of function fields, particularly in relation to recent developments in anabelian
geometry (see \cite{BogTsc12}). One of the spectacular achievements of the latter
is that (see \cite[Theorem 13.2]{Bog91}) abelian subgroups of absolute Galois
groups $G = \Gal(K^s/K)$ are determined already by considering the Galois group only up
to the second term of its lower central series,
$G^c = G / [[G,G],G]$.
Passing to a maximal pro-$p$ quotient $\mathcal{G}_p$ of $G$, we find
the $p$-part of the Bogomolov multiplier as
\[
\B_0(G)_p = \B_0(\mathcal{G}_p) = 
\B_0(\mathcal{G}_p^c) = 
\varinjlim_{N \unlhd_o \mathcal{G}_p^c} \B_0(\mathcal{G}_p^c / N).
\]
Understanding these Bogomolov multipliers is therefore reduced to
considering either pro-$p$ groups of nilpotency class $2$, or equivalently
their finite quotients. 

\subsection{Contributions and reader's guide}

Motivated by the above, 
the aim of this paper is to investigate {\em generic behaviour of Bogomolov
multipliers} of nilpotent pro-$p$ groups of sufficiently small nilpotency
class, especially those of nilpotency class $2$, as well as their finite
quotients. In doing so, we are able to
extract the following rather surprising result (Theorem~\ref{t:nontrivial_bog_1}).

\begin{theorem*}[$n \to \infty$]

Fix a prime $p > 2$ and a positive integer $M$.
Let $\#_{all}(n)$ be the number of all $p$-groups of order $p^n$, and let
$\#_{\B_0 \geq M}(n)$ be the number of those $p$-groups of order $p^n$ whose Bogomolov
multiplier is of order at least $M$. Then
\[
\lim_{n \to \infty}
\frac
{\log_p \#_{\B_0 \geq M}(n)}
{\log_p \#_{all}(n)}
= 1.
\]
\end{theorem*}

\begin{corollary*}
Quotient varieties of faithful representations of 
\mbox{log-generic $p$-groups} are not stably rational.
\end{corollary*}

This is quite in contrast with thinking that groups 
possessing the Bogomolov $2$-cocycles are rather special.
We now explain our technique together with side applications.

\smallskip

In Section \ref{section:overture} we first
dualize the Bogomolov multiplier $\B_0(G)$ into its homological
version (Subsection \ref{ss:homological_intepretation}). 
This dualized version appears naturally in non-commutative 
Iwasawa theory as $\SK_1(\ZZ_p \llbracket G \rrbracket)$
(Subsection \ref{ss:iwasawa_theory}).
We also give some remarks on general structural properties of this object
(Subsection \ref{ss:rank_and_exponent}).

\smallskip

In Section \ref{section:via_lie} we 
restrict to nilpotent pro-$p$ groups with a special regard to
Lie groups over the $p$-adic field $\QQ_p$.
We transport the dualized object from the previous section to the
category of Lie algebras (Subsections \ref{ss:lazard_correspondence}
and \ref{ss:transporting_sk1})
by exploiting Hopf formulae for Schur multipliers
(Subsection \ref{ss:hopf_formulae}).
We then restrict to objects of nilpotency class at most $2$
and no torsion (Subsection \ref{ss:class_2}).
We show how to parametrize such objects by using the 
Grassmannian variety  $\Gr(r,\binom{d}{2})$
as the moduli space (Subsection \ref{ss:parametrizing_lie_algebras}). 
To every point $L \in \Gr(r, \binom{d}{2})$ we associate a
Lie algebra $\Lie$ and then a group $G_L$ via the Lazard correspondence.
Here, $d$ and $r$ are the natural parameters encoding
the number of generators and essential relators
of the objects.
Using the above, we establish a Lie criterion
over $\QQ_p$ for finiteness of the Bogomolov multiplier
(Subsection \ref{ss:sk1_of_lie_algebra}).

\smallskip

In Section \ref{section:decomposability_Qp} we 
express the general Bogomolov multiplier $\B_0(G_L)$
in terms of the image of some rational map that we call the decomposability
map (Subsection~\ref{ss:decomposability_morphism_Qp}). 
It is defined on some quasiprojective variety over $\QQ_p$ 
and maps into the moduli space. 
This map is the heart of this paper.
We find a subvariety of the
moduli space of small codimension (more precisely, the ratio of its codimension and the dimension of the moduli space is asymptotically zero) that misses the image of this
map and thus consists entirely of
groups with large Bogomolov multipliers
(Subsection \ref{ss:plenty_of_nondecomposable_subspaces}). 
We study the generic
behaviour by analysing the local behaviour of this map using methods of
differential geometry (Subsection \ref{ss:differential_of_decomposability}), and we give a precise criterion on when the tangent map
is generically injective or surjective in terms of $r$ and $d$
(Subsections \ref{ss:submersiveness} and \ref{ss:immersiveness}).
The crucial step here is analysing the case $r = \binom{d-2}{2} + 1$,
when the decomposability map turns out to be a local
isomorphism.

\smallskip

In Section \ref{section:decomposability_Fp} we finally go back to the original
motivation and transfer results from the analytic world
to their finite quotients of exponent $p$ (Subsections \ref{ss:groups_and_grassmannian}
and \ref{ss:decomposability_morphism_Fp}).
Groups with few
generators are dealt with separately, where we can exploit the presence of
an action of a general linear group due to the dimensions being sufficiently 
small (Subsection \ref{ss:few_generators}).
Some of the general phenomena are visible in these calculations,
particularly when the decomposability map is a local isomorphism.
In order to deal with general groups, we expand the
map describing the general Bogomolov multiplier into a scheme map
(Subsection \ref{ss:asymptotics_of_decomposability}).
The value of the parameters $r$ and $d$ that is most relevant for
our understanding of generic $p$-groups corresponds to
the situation when the $p$-adic base change of the 
decomposability map is a local surjection.
We use some arithmetic geometry
to then conclude that in this case,
the proportion of elements in the image of the relevant
map over the finite field $\FF_p$ is bounded away from $0$ as $p$
tends to infinity (Theorem \ref{t:large_r_asymptotics}).

\begin{theorem*}[$p \to \infty$]

Fix $r$ and $d$ with $r \geq \binom{d-2}{2} + 1$.
Then
\[
\liminf_{p \to \infty} \frac{|\{ L \in \Gr(r, \binom{d}{2}) \mid \B_0(G_L) = 0 \}|}{|\Gr(r, \binom{d}{2})|} \geq
\left( \frac{1}{C_{d-2}} \right) ^r,
\]
where $C_{d-2}$ is the Catalan number.

\end{theorem*}

Thus there are many groups with trivial Bogomolov multipliers. The first issue with
the above result is that it only applies to large enough primes $p$, while
we are interested in finite $p$-groups for a fixed $p$. 
The second issue is that even if the above inequality would hold at every fixed
$p$, it would still not enable us to provide numerically many
finite $p$-groups $G_L$ with vanishing Bogomolov multipliers. The reason
for this is that the constant on the RHS of the above theorem converges to
zero too quickly. 
In fact, it comes as a surprise that just the opposite can be achieved (Subsection \ref{ss:log_generic_groups}).
By fixing $p$ and letting the parameters $r$ and $d$ suitably tend to infinity, 
the produced finite groups $G_L$ will be of order $p^n$ with $n$ depending only
on $r$ and $d$.
The parameters can be set in such a way that the number of $p$-groups
produced is log-comparable to the number of {\em all} groups of order $p^n$,
and we show that with increasing $n$, log-generic groups $G_L$ have
non-trivial Bogomolov multipliers (Theorem $(n \to \infty)$ / \ref{t:nontrivial_bog_1}). 
This relies on transferring the
small codimension subvariety mentioned above from the analytic world into finite fields. Whence log-generic $p$-groups produce quotient varieties that are very much not stably rational.

\smallskip

The paper concludes with an application of our technique to another problem
involving commutators (Subsection \ref{ss:commutators_in_the_end}). We show
that in log-generic  $p$-groups, not every element of the derived subgroup is
a simple commutator.

\subsection*{Acknowledgements}

The authors thank 
Oihana Garaialde Ocaña for conversations on the Lazard correspondence and
many other nice remarks,
Jon González Sánchez for providing the preprint \cite{GSWeixx},
Immanuel Halupczok for sharing his knowledge of general algebraic geometry,
Geoffrey Janssens for helpful thoughts on the exposition,
and
Daniel Loughran for pointing out the relevant version of the fibre dimension theorem and patiently responding to related questions.

\section{Overture}
\label{section:overture}

\noindent This section is devoted to recognizing the dual of the 
Bogomolov multiplier within Iwasawa theory and collecting some known
and unknown structural properties of this object. We will henceforth 
switch depending on the context between this homological interpretation 
and the original Bogomolov multiplier.

\subsection{Iwasawa theory} \label{ss:iwasawa_theory}

Let $G$ be a compact analytic pro-$p$ group over $\QQ_p$,
and assume $p > 2$. 
Associated to $G$ is its completed group ring
\[
\ZZ_p \llbracket G \rrbracket =
\varprojlim_{N \lhd_o G} \ZZ_p[G/N]
\]
over the $p$-adic integers $\ZZ_p$, where the inverse limit runs over all open normal subgroups $N$ of $G$. {\em Non-commutative Iwasawa theory} deals with
investigating this ring from the point of view of {\em $\K$-theory}. More precisely,
one is interested in understanding the abelian group $\K_1(\ZZ_p \llbracket G
\rrbracket )$ in relation with constructing $L$-functions for $p$-adic
representations (see \cite{Wit17}). The continuous representation theory of 
$G$ is captured
by the completed group algebra over the $p$-adic field
\[
\QQ_p \llbracket G \rrbracket =
\varprojlim_{N \lhd_o G} \QQ_p[G/N].
\]
This algebra decomposes into a product of matrix algebras over finite
extensions of $\QQ_p$ (see \cite[Section 3]{SchVen13split}), making it easy
to understand its $\K_1$-group. There is a natural scalar extension map
\[
\K_1(\ZZ_p \llbracket G \rrbracket)
\to
\K_1(\QQ_p \llbracket G \rrbracket).
\]
The deficit of this map being injective is measured by its kernel
$\SK_1(\ZZ_p \llbracket G \rrbracket)$. 
This object can be represented as a limit
\[
\SK_1(\ZZ_p \llbracket G \rrbracket) =
\varprojlim_{N \lhd_o G} \SK_1(\ZZ_p[G/N]).
\]
It follows from \cite[Theorem 3]{Oli80} that the groups $\SK_1(\ZZ_p[G/N])$ are
finite $p$-groups, and so $\SK_1(\ZZ_p \llbracket G \rrbracket)$ is an abelian
pro-$p$ group. 

\subsection{Dualizing into the Bogomolov multiplier}

There is a homological way of interpreting the finite
$\SK_1$'s, studied in detail by Oliver via
the $p$-adic logarithm (see \cite{Oli80}).
This description can be dualized (see \cite{Mor12})
and extended to continuous Galois cohomology (see 
\cite[Corollary 2.2]{SchVen13lie}), 
yielding the following isomorphism:
\begin{equation} \label{e:sk1_as_kernel}
\SK_1(\ZZ_p \llbracket G \rrbracket)^* \cong
\bigcap_{A \leq G, \, [A,A] = 1} \ker \left( 
\HH^2(G, \QQ_p/\ZZ_p)
\xrightarrow{\text{$res^G_A$}}
\HH^2(A, \QQ_p/\ZZ_p)
\right).
\end{equation}
Here the dual is taken with respect to $\QQ_p/\ZZ_p$.
Thus $\SK_1(\ZZ_p \llbracket G \rrbracket)^*$ is nothing
but $\B_0(G)$.

\subsection{Homological interpretation} \label{ss:homological_intepretation}

We will work with the homological version of $\B_0(G)$.
This can be obtained directly by dualizing
\eqref{e:sk1_as_kernel}.  The abelian pro-$p$ group  $\SK_1(\ZZ_p \llbracket G
\rrbracket)$ coincides with its double dual (see \cite[Theorem
2.9.6]{RibZal00}). Thus we obtain an exact sequence
\[
\left( \prod_{A \leq G, [A,A] = 1}  \HH^2(A, \QQ_p/\ZZ_p) \right)^*
\xrightarrow{\text{$\left( \prod res^G_A \right)^*$}}
\HH^2(G, \QQ_p/\ZZ_p)^*
\xrightarrow{}
\SK_1(\ZZ_p \llbracket G \rrbracket)
\xrightarrow{}
0.
\]
The Universal Coefficient Theorem (see \cite[Section 7.6]{Por09}) gives a natural
isomorphism $\HH^2(G,\QQ_p/\ZZ_p)^* \cong \HH_2(G,\ZZ_p)$,
and therefore we have an exact sequence 
\begin{equation} \label{e:sk1_as_quotient}
\bigoplus_{A \leq G, [A,A] = 1}  \HH_2(A, \ZZ_p)
\xrightarrow{\text{$\bigoplus cores^G_A$}}
\HH_2(G, \ZZ_p)
\xrightarrow{}
\SK_1(\ZZ_p \llbracket G \rrbracket)
\xrightarrow{}
0.
\end{equation}

\subsection{Exterior squares and commutator relations}

As in the case of finite groups (see \cite[Theorem 3.2]{Mor12}), the sequence
\eqref{e:sk1_as_quotient}
can be interpreted via the notion of {\em non-abelian exterior squares} (see
\cite[Section 8.6]{Por09}). The description goes as follows. 
Let $G \pwedge G$ be the profinite group topologically generated by symbols
$g \wedge h$ for $g, h \in G$ subject to the universal relations
\[
xy\wedge z = (x^y\wedge z^y)(y\wedge z), \quad 
x\wedge yz = (x\wedge z)(x^z\wedge y^z), \quad
x\wedge x = 1
\]
for all $x, y, z \in G$, where $x^y$ denotes the conjugate $y^{-1} x y$.
This group may be thought of as the universal model
for the derived subgroup $[G,G]$. Correspondingly, there is a commutator map
\begin{equation} \label{e:commutator_map_group}
\kappa \colon G \pwedge G \to [G,G], \quad g \pwedge h \mapsto [g,h].
\end{equation}
The defect that the above relations between wedges are in fact all relations
between commutators is measured by the kernel $\M(G) := \ker \kappa$.
It follows from \cite[Section 8.6.5]{Por09} that there is a natural isomorphism
$\M(G) \cong \HH_2(G, \ZZ_p)$. Moreover, corestrictions from abelian subgroups
$A \leq G$ correspond to natural maps $A \pwedge A \to G \pwedge G$.
Set $\M_0(G) := \langle g \wedge h \mid g, h \in G, \, [g,h] = 1 \rangle$.
Equation \eqref{e:sk1_as_quotient} can therefore be interpreted as
\begin{equation} \label{e:sk1_as_kernel_of_commutator_map_group}
\SK_1(\ZZ_p \llbracket G \rrbracket) \cong \frac{\M(G)}{\M_0(G)},
\end{equation}
the group of relations between commutators in $G$ modulo the universal relations
and commuting pairs (see \cite{JezMor15unicp}).

\subsection{Bounding the rank and exponent} \label{ss:rank_and_exponent}

Both for the purposes of Iwasawa theory and the rationality problem, it is of
particular interest to know whether or not the object $\SK_1(\ZZ_p \llbracket G
\rrbracket)$ is trivial.  Exploiting the fact that in our context, the group $G$
is of {\em finite rank} (this means that there is an upper bound on the 
minimal number of generators of its closed subgroups, see \cite[Corollary 8.34]{DDMS99}), we provide some general
bounds on the structure of this object. Recall that such a group $G$ 
is equivalently described as a finitely generated pro-$p$ group that possesses
a finite index subgroup $H$ that is {\em uniform}. The latter means that
$H$ is finitely generated, torsion-free and powerful, i.e., $[H,H] \leq H^{p}$
(see \cite{DDMS99}).

We first show that the rank of the Bogomolov multiplier can be bounded in terms of the rank of $G$.

\begin{proposition} \label{p:sk1_fingen}
$\SK_1(\ZZ_p \llbracket G \rrbracket)$ is a finitely generated abelian pro-$p$ group.
\end{proposition}
\proof
Let $r$ be the rank of the group $G$.
All finite quotients $G/N$ as well as
their subgroups can be generated by $r$ elements. Furthermore, every such 
finite quotient contains a powerful subgroup of index bounded by a function of
$r$ independently of $N$ (see \cite[Theorem 2.13]{DDMS99}). 
The rank of the second integral homology of this powerful subgroup can be bounded
in terms of $r$ (see \cite[Corollary 2.2]{LubMan87}),
and therefore the same is true for its Bogomolov multiplier.
Passing from this subgroup to $G$, the rank of the multiplier gets at worst multiplied by a constant depending only on its
index and therefore only on $r$ (see \cite[Proposition 6.1]{JezMor15cpext}).
Thus the rank of the Bogomolov multiplier of the finite quotients $G/N$ is
bounded by a function of $r$. It follows that the limit
$\SK_1(\ZZ_p \llbracket G \rrbracket)$ is topologically finitely
generated (see \cite[Proposition 4.2.1]{Wil98}).
\endproof

A similar argument proves that the exponents of Bogomolov multipliers of
finite quotients of $p$-adic analytic groups can asymptotically be controlled.

\begin{proposition}
There exists a function $f(r)$ such that for every pro-$p$ $p$-adic
analytic group $G$ of rank $r$, we have
$\exp \B_0(G/N) \leq f(r) \cdot \exp G/N$ for every $N \lhd_o G$.
\end{proposition}
\proof
As in the proof of Proposition \ref{p:sk1_fingen}, there is a powerful subgroup of $G/N$ of rank and index
bounded in terms of $r$. The exponent of the Bogomolov multiplier of such a 
subgroup is bounded by the exponent of $G/N$ (see \cite[Theorem 2.4]{LubMan87}). It follows that
the exponent of $\B_0(G/N)$ can be bounded by a function of $r$ multiplied by
$\exp G/N$ (see \cite[Proposition 6.2]{JezMor15cpext}).
\endproof

It follows from the above proofs that when the uniform subgroup of $G$ is
abelian, the group $\SK_1(\ZZ_p \llbracket G \rrbracket)$ is of finite exponent
bounded in terms of the rank of $G$, and so it is a finite $p$-group. This is
the case, for example, when $G$ is a pro-$p$ group of finite coclass (see
\cite[Theorem 10.1]{DDMS99} and \cite[Theorem 4.4]{FerJez16}).

\section{Via Lie Theory}
\label{section:via_lie}


\noindent In this section we linearise the description of the Bogomolov
multiplier to Lie algebras of uniform pro-$p$ groups.  This description,
however, is not entirely equivalent in the categorical sense to the Bogomolov
multiplier. In order to properly transfer homology, we will need to work with
nilpotent pro-$p$ groups of bounded class. Relying on Hopf formulae for both
groups and Lie algebras, it is possible to transfer the whole Schur multiplier
and thus also the Bogomolov multiplier. For nilpotent uniform pro-$p$ groups,
we obtain a simple Lie criterion for triviality of the Bogomolov multiplier,
that does coincide with the naive notion introduced at the beginning of this
section.

\subsection{The Lie algebra and the commutator map}

Let $G$ be a uniform pro-$p$ group. One can associate to it a Lie
algebra $\Lie = \log G$ over $\ZZ_p$ (see \cite[Chapter 9]{DDMS99}).
The Lie bracket $[\;,\;]_L$ determines a mapping
\[
\kappa_L \colon \Lie \wedge \Lie \to \Lie, \quad x \wedge y \mapsto [x,y]_L.
\]
This map may be thought of as a linearised version of \eqref{e:commutator_map_group}.

\subsection{Decomposable wedges and $\SK_1$}

One can linearise the description
\eqref{e:sk1_as_kernel_of_commutator_map_group} of $\SK_1(\ZZ_p \llbracket G
\rrbracket)$. Let 
\[
D\Lie = \{ x \wedge y \mid x,y \in \Lie \}
\subseteq \Lie \wedge \Lie 
\]
denote the set of
{\em decomposable wedges} in $\Lie$.  It is clear that whenever
$x,y \in \Lie$ commute, we have $x \wedge y \in D \Lie \cap \ker \kappa_L$.

It is shown in \cite[Theorem 5.5]{SchVen13lie} that $\SK_1(\ZZ_p
\llbracket G \rrbracket)$ is trivial if and only if $\ker \kappa_L$ can be
generated by elements of $D \Lie \cap \ker \kappa_L$. This is done by passing to
$\FF_p$-coefficients in cohomology and using quite heavy calculations relating
the structure of $G$ to that of $\Lie$. Thus the intuitive way of transferring 
the Bogomolov multiplier to Lie algebras makes sense. However, this only
works with the triviality criterion. We will now transfer the whole object
and thus provide a more conceptual clarification of this triviality criterion.

\subsection{The Lazard correspondence} \label{ss:lazard_correspondence}

We can avoid the above mentioned calculations by restricting
to the class of nilpotent pro-$p$ groups of small nilpotency class by
using the {\em Lazard correspondence} (see \cite{Laz54}). 
This is an isomorphism of categories
between the category of finitely generated nilpotent pro-$p$ groups of
nilpotency class less than $p$ and finitely generated nilpotent Lie algebras
over $\mathbb{Z}_p$ of nilpotency class less than $p$. 
To a Lie algebra $\Lie$, one associates the pro-$p$ group $\exp \Lie$,
which coincides with $\Lie$ as a set and whose multiplication arises from the Baker-Campbell-Hausdorff formula
\[
x \cdot y = \log ( \exp(x)  \cdot \exp(y))
\]
for $x, y \in \Lie$. The logarithm and exponential here arise from formal power
series, truncated from $p$ on due to the restriction on the nilpotency class.
Conversely, to a pro-$p$ group $G$ one associates the Lie algebra $\log G$ 
which coincides with $G$ as a set and with
\[
x + y = \exp ( \log(x) + \log(y) ),
\quad
[x,y] = \exp ( [\log(x), \log(y)] )
\]
for $x, y \in G$. The functors $\exp$ and $\log$ are inverse to each other
and induce an isomorphism of the relevant categories. They transform subgroups
to subalgebras, normal subgroups to ideals, commuting pairs to
commuting pairs, centre to centre, generating
sets to generating sets, etc.
We remark that these functors are the truncated versions of functors that transform
uniform groups to their Lie algebras and vice versa.

\subsection{Hopf formulae} \label{ss:hopf_formulae}

Our objective is to transport the functor $\SK_1$ from the category of groups
to that of Lie algebras. In light of \eqref{e:sk1_as_quotient}, we will do
this by first transporting the second homology. This relies on the following
universal expression of the second homology in both group and Lie algebra
categories.

\begin{lemma}[Hopf formulae]
Let $F / R$ be a free presentation of either a pro-$p$ group $G$
or a $\mathbb{Z}_p$-Lie algebra $\Lie$. Then
\[
\HH_2(F/R, \ZZ_p) \cong \frac{[F,F] \cap R}{[F,R]}.
\]
\end{lemma}
\proof
This is classical for groups. See 
\cite[Section 8.2.3]{Por09} for pro-$p$ groups,
the heart of the argument being the $5$-term exact sequence
from the LHS spectral sequence. 
The same works for Lie algebras, see also \cite{Ell87}.
A universal explanation of this phenomenon is given in \cite{BDIL07}.
\endproof

\subsection{Universal central extensions}

The objects described by Hopf formulae can be interpreted as universal
objects associated to the group or the Lie algebra in the following sense.
Let $F/R$ be a free presentation of a finitely generated pro-$p$ group $G$
with $\dd(F) = \dd(G)$. 
Here $\dd(G)$ is the minimal number of topological generators of $G$.
There exists a largest $\dd(G)$-generated central
extension of $G$, and this is precisely the extension
\[
1 \to \frac{R}{[F,R]} \to \frac{F}{[F,R]} \to G \to 1.
\]

\subsection{Transporting homology}
Denote $E = F / [F,R]$ and let $p \colon E \to G$ be the projection.
The group $E$ comes equipped with the additional projection 
$a \colon E \to E / [E,E]$. Combine the two projections $p, a$
into a product map 
\[
\frac{F}{[F,R]} = E \to G \times (E / [E,E]) = \frac{F}{R} \times \frac{F}{[F,F]}
\]
and let $K$ be the kernel
of this map. Note that $K = ([F,F] \cap R)/[F,R]$.
Therefore the object $K \cong \HH_2(G, \mathbb{Z}_p)$ can be expressed in purely
categorical terms. Note that since we are assuming $G$ to be finitely generated,
all the relevant objects of the construction stay in the category of finitely
generated nilpotent pro-$p$ groups. Moreover, if $G$ is of nilpotency class $c$,
then $E$ is of nilpotency class at most $c + 1$. Using the Lazard isomorphism
of categories, we can therefore fully transport the construction, starting from
the universal central extension. In order for the correspondence to apply to
all the relevant objects, we need to assume that $c + 1 < p$.
In particular, starting with a group of nilpotency class $2$, this applies
to all primes $p > 3$.
Finally we obtain the natural isomorphism 
\[
\log \HH_2(G, \ZZ_p) \cong \HH_2(\Lie, \ZZ_p).
\]

\begin{remark}
Let $G$ be a uniform pro-$p$ group, not necessarily nilpotent.
In this situation, one can still define a Lazard correspondence.
The problem of making our construction work in this more general
setting is that the additional object $E$ that appears in the construction
might not be uniform and it is not clear what its associated
Lie algebra should be. For a unified explanation of why the above homological
transport works, one would therefore need a more general Lazard correspondence.
See \cite[Section 9.4]{GSWeixx}.
\end{remark}

\begin{remark}
The above argument of transporting homology
does not work for finite $p$-groups, since the object $E$ falls out
of the relevant category. It is nevertheless possible to obtain
the same result for the case of finite $p$-groups,
either using Schur's theory of covers (see \cite{EHZ12})
or more explicitly via the usual interpretations of low-dimensional
cohomology (see \cite{GarGon17}).
\end{remark}


\subsection{Transporting $\SK_1$} \label{ss:transporting_sk1}

Given an abelian subgroup $A \leq G$, there is a corestriction
map $cores^G_A \colon \HH_2(A, \ZZ_p) \to \HH_2(G, \ZZ_p)$. It follows from the above
that there is an induced
corestriction map
$\log(cores^G_A) \colon \HH_2(\log A, \ZZ_p) \to \HH_2(\log G, \ZZ_p)$.
It follows from \eqref{e:sk1_as_quotient}
that
\[
\SK_1(\ZZ_p \llbracket G \rrbracket) \cong \frac{\HH_2(G, \ZZ_p)}
{\langle \im cores^G_A \mid A \leq G, \, [A,A] = 1 \rangle},
\]
and the latter can be transferred with $\log$ to
\[
\log(\SK_1(\ZZ_p \llbracket G \rrbracket)) \cong \frac{\HH_2(\Lie, \ZZ_p)}
{\langle \im \log(cores^G_A) \mid A \leq G, \, [A,A] = 1 \rangle}.
\]
Abelian subgroups of $G$ correspond precisely to abelian Lie subalgebras
of $\Lie$, and so we obtain
\[
\log(\SK_1(\ZZ_p \llbracket G \rrbracket)) \cong \frac{\HH_2(\Lie, \ZZ_p)}
{\langle \im cores^{\Lie}_{\mathcal{A}} \mid \mathcal{A} \leq \Lie, \, 
[\mathcal{A},\mathcal{A}] = 1 \rangle}.
\]
Following \cite[Remark 9.19]{GSWeixx}, there is a direct way of computing
$\HH_2(\Lie, \ZZ_p)$ via the standard complex
(see \cite[Theorem 4.2]{HiltonStammbach}). To do this, 
consider the Jacobi map
\[
\theta_{\Lie} \colon \Lie \wedge \Lie \wedge \Lie \to \Lie \wedge \Lie, \qquad
x \wedge y \wedge z \mapsto [x, y] \wedge z + [y, z] \wedge x + [z, x] \wedge y.
\]
We have a natural isomorphism
\[
\HH_2(\Lie, \ZZ_p) \cong \frac{\ker \kappa_{\Lie}}{\im \theta_{\Lie}}
\]
For an abelian Lie subalgebra $\mathcal{A} \leq \Lie$, we therefore
have $\HH_2(\mathcal{A}, \ZZ_p) \cong \ker \kappa_{\mathcal{A}} =
\mathcal{A} \wedge \mathcal{A}$. We thus obtain a natural isomorphism
\begin{equation} \label{e:sk1_definition_lie_algebra}
\log(\SK_1(\ZZ_p \llbracket G \rrbracket)) \cong \frac{\ker \kappa_{\Lie}}
{\im \theta_{\Lie} + \langle D \Lie \cap \ker \kappa_{\Lie} \rangle }
\end{equation}
and we correspondingly denote the left hand side by $\SK_1(\Lie)$.
We have therefore transferred $\SK_1$ into the category of nilpotent 
Lie algebras over $\ZZ_p$ of nilpotency class less than $p-2$.

\begin{remark}
The object $\SK_1(\Lie)$ as we have defined it makes sense for any Lie algebra
$\Lie$. In fact, it coincides with the quotient 
(in the sense of \eqref{e:sk1_as_kernel_of_commutator_map_group})
of the Lie bracket induced map $\ker(\Lie \wedge \Lie \to \Lie)$,
where $\wedge$ is the non-abelian exterior square of the Lie algebra
$\Lie$ (see \cite{Ell87}).
This object has been explored independently in \cite{RPN18} (see also \cite{RPN20} and \cite{CheMa21}).
\end{remark}

\subsection{Triviality of $\SK_1$}

We show how the above can be utilized in order to give a quick argument that
triviality of $\SK_1(\ZZ_p \llbracket G \rrbracket)$ is precisely detected
by decomposability of wedges (see \cite[Theorem 5.5]{SchVen13lie}).
We are assuming that $G$ is both a uniform pro-$p$ group and that
$G$ is nilpotent of nilpotency class at most $p-2$.

Note that since $\Lie$ is a powerful Lie algebra, we have
$[\Lie, \Lie] \subseteq p \Lie$, and therefore
$\im \theta_{\Lie} \subseteq p \Lie \wedge \Lie \cap \ker \kappa_{\Lie}$.
The algebra $[\Lie, \Lie] \cong \Lie \wedge \Lie / \ker \kappa_{\Lie}$
is (additively) torsion-free, and so $\im \theta_{\Lie} \subseteq p \ker \kappa_{\Lie}$.
Note that $\SK_1(\Lie)$ is trivial if and only $\ker \kappa_{\Lie}
= \im \theta_{\Lie} + \langle D \Lie \cap \ker \kappa_{\Lie} \rangle$.
Since $\im \theta_{\Lie}$ is contained in the Frattini subgroup of $\ker \kappa_{\Lie}$,
we have the following.

\begin{corollary}
Let $G$ be a uniform pro-$p$ group of nilpotency class at most $p-2$.
Then
$\SK_1(\ZZ_p \llbracket G \rrbracket) = 0$ if and only if
$\ker \kappa_{\Lie} = \langle D \Lie \cap \ker \kappa_{\Lie} \rangle$.
\end{corollary}

\subsection{Nilpotency class $2$} \label{ss:class_2}

We now specialize all of the above to groups and algebras of nilpotency class
$2$. Thus we assume throughout that $p > 3$ and we restrict to the category of
pro-$p$ groups that are finitely generated and
nilpotent of class at most $2$. We will additionally assume that our groups
and algebras are torsion-free in order to pass from $\ZZ_p$ to $\QQ_p$,
and we will also assume that their abelianizations are torsion-free
in order to have a uniform description.
On the other hand, our construction will not need the more powerful restriction that
the groups and algebras are powerful.
Our final objective will be to analyse the behaviour of $\SK_1$ in the generic case.
More precisely, we want to consider the Bogomolov multiplier of 
a random pro-$p$ group from the above category.
The main reason why we restrict to objects of nilpotency class $2$
is that these can be parametrized in a natural way as we explain below,
and this makes it possible to talk about random objects.

\subsection{Parametrizing Lie algebras of nilpotency class $2$}
\label{ss:parametrizing_lie_algebras}

Let $V$ be a finitely-generated free $\ZZ_p$-module. Consider $W = V \wedge V$ and let
$L$ be a $\ZZ_p$-submodule of $W$. The natural projection map
$\pi \colon W \to W/L$ can be thought of as determining the Lie product
on the $\ZZ_p$-module $\mathcal{L} = V \oplus W/L$,
\[
[v_1 + w_1, v_2 + w_2]_{\Lie} = \pi(v_1 \wedge v_2).
\]
Thus we obtain a Lie algebra of nilpotency class at most $2$. 
Conversely, every Lie algebra $\Lie$
of nilpotency class $2$
with $\Lie / [\Lie, \Lie]$ torsion-free is obtained in this way,
where we take $L = \ker (\Lie \wedge \Lie \to [\Lie, \Lie])$
to be the kernel of the Lie bracket map.

We will require the algebra $\mathcal{L}$ to be torsion-free.
This occurs precisely when $W/L$ has no torsion.
We will call such $\ZZ_p$-submodules $L \leq W$ co-torsion-free.
The following lemma shows that such submodules naturally correspond to
$\QQ_p$-subspaces of $W \otimes \QQ_p$.
We will prefer to work in this setting.

\begin{lemma} \label{l:bijection_subspaces_submodules} 
There is a bijection between $r$-dimensional vector subspaces of $\QQ_p^k$ and
co-torsion-free $\ZZ_p$-submodules of $\ZZ_p^k$ of rank $r$.
The bijection is given as
\[
\QQ_p^k \supseteq L \mapsto L \cap \ZZ_p^k, \qquad \ZZ_p^k \supseteq L \mapsto L \otimes \QQ_p.
\]
\end{lemma}
\proof
Let us first verify that the map is well-defined. Given $L \subseteq \QQ_p^k$,
the $\ZZ_p$-module $L \cap \ZZ_p^k$ is a  torsion-free abelian pro-$p$ group.
Its rank equals the rank of $(L \cap \ZZ_p^k) \otimes \QQ_p$. Since there is an
$n \geq 0$ such that $p^n L \subseteq \ZZ_p^k$,  the latter rank is equal to the
dimension of $L$, that is $r$. Let us now verify that $\ZZ_p^k / (L \cap
\ZZ_p^k)$ is torsion-free. Were there an element $\ell \in \ZZ_p^k$ with $p \ell
\in L \cap \ZZ_p^k$, then $p \ell \in L$ and so $\ell \in L$, and so $\ell \in L
\cap \ZZ_p^k$.

We now check that the above map is a bijection. To this end, we only need
to verify that for a co-torsion-free $\ZZ_p$-module $L \leq \ZZ_p^k$,
we have $(L \otimes \QQ_p) \cap \ZZ_p^k = L$. Let $\ell = \sum_i \alpha_i \ell_i
\in (L \otimes \QQ_p) \cap \ZZ_p^k$ with $\alpha_i \in \QQ_p$ and $\ell_i \in L$.
Therefore there is an $n \geq 0$ such that $p^n \alpha_i \in \ZZ_p$ for all $i$,
and so $p^n \ell \in \sum_i \ZZ_p \ell_i \subseteq L$. It follows from the
assumption on $L$ that $\ell \in L$.
\endproof

Once such a subspace $L$ is chosen, we can construct the corresponding Lie
algebra $\mathcal{L}$.  In this setting, it is also easy to produce  uniform
algebras or groups.  This is achieved by simply replacing an algebra
$\mathcal{L}$ as above with $p \mathcal{L}$.

All in all, we have designed a map that takes as input a subspace
$L$ of dimension $r$ in the fixed $\QQ_p$-vector space $W$
and constructs a Lie algebra of nilpotency class
$2$ by imposing elements of $L \cap \ZZ_p^k$ as relators between Lie brackets.
The subspaces $L$ can be thought of as belonging to the Grassmannian variety
$\Gr(r, W)$ (see \cite{Grassmannian}). Thus we have a parametrization
from $\Gr(r, W)$ to the set $\mathbf{N_2Lie}$ of  Lie algebras of
nilpotency class $2$ that are torsion-free and have torsion-free abelianizations:
\begin{equation} \label{eq:parametrization_Qp}
\grasslie \colon \Gr(r, W) \to \mathbf{N_2Lie}, \qquad L \mapsto \Lie.
\end{equation}
To every such a Lie algebra $\Lie$ we associate via the Lazard correspondence
a group that we denote $\exp \Lie = G_L$. 
The map $\grasslie$ is therefore also parametrizing
the relevant groups.

\subsection{$\SK_1$ of the Lie algebra} \label{ss:sk1_of_lie_algebra}

Given a subspace $L \subseteq W$, we determine the $\SK_1$ of the Lie algebra
$\Lie = \grasslie(L)$ using
\eqref{e:sk1_definition_lie_algebra}. Note that in this case,
$\im \theta_{\Lie}
\subseteq \langle D \Lie \cap \ker \kappa_{\Lie} \rangle$,
and we therefore have
\[
\SK_1(\Lie) =
\frac{L}{\langle D \Lie \cap L \rangle}.
\]
Nothing changes if we replace the Lie algebra by the powerful one $p \Lie$,
\[
\SK_1(p \Lie) =
\frac{\ker \kappa_{p  \Lie}}
{\langle D(p \Lie) \cap \ker \kappa_{p  \Lie} \rangle} =
\frac{p^2 \ker \kappa_{ \Lie} }
{\langle p^2 \left( D  \Lie \cap \ker \kappa_{ \Lie} \right) \rangle} \cong
\frac{L}{\langle D \Lie \cap L \rangle}.
\]
Using the correspondence from Lemma \ref{l:bijection_subspaces_submodules},
we can return from submodules to subspaces and consider the Lie algebra
$\Lie \otimes \QQ_p$ with its own Bogomolov multiplier
\[
\SK_1(\Lie \otimes \QQ_p) = 
\frac{L \otimes \QQ_p}{\langle D (\Lie \otimes \QQ_p) \cap (L \otimes \QQ_p) \rangle}.
\]
The following lemma establishes the connection between the two objects.

\begin{lemma} \label{l:sk1_subspaces_submodules}
$\SK_1(\Lie)$ is finite if and only if $\SK_1(\Lie \otimes \QQ_p) = 0$.
\end{lemma}
\proof
Suppose first that $\SK_1(\Lie)$ is finite. Thus there exists an $n \geq 0$
such that for every element $\ell$ of the $\ZZ_p$-basis of $L$, we have $p^n
\ell \in \langle D(\Lie) \cap L \rangle$. It follows that $\ell \in \langle
D(\Lie \otimes \QQ_p) \cap (L \otimes \QQ_p) \rangle$. This implies that $L
\subseteq \langle D(\Lie \otimes \QQ_p) \cap (L \otimes \QQ_p) \rangle$, and so
$\SK_1(\Lie \otimes \QQ_p) = 0$.

Conversely, assume that $\SK_1(\Lie \otimes \QQ_p) = 0$. Then every element
$\ell$ of the $\ZZ_p$-basis of $L$ belongs to  $\langle D(\Lie \otimes \QQ_p)
\cap (L \otimes \QQ_p) \rangle$. Thus there is an $n \geq 0$ such that $p^n \ell
\in \langle D(\Lie) \cap L \rangle$ for all basis elements of $L$, and
therefore $p^n L \subseteq \langle D(\Lie) \cap L \rangle$. This means that
$\SK_1(\Lie)$ is a finitely generated torsion abelian pro-$p$ group, and so
finite.
\endproof

Here is a neat example showing how torsion can indeed occur.

\begin{example} \label{ex:torsion_in_sk1}
Consider the case when $V$ is the free $\ZZ_p$-module of dimension $4$
generated by $\{ e_1, e_2, e_3, e_4 \}$. Let $L$ be the $\ZZ_p$ submodule
\[
L = \langle 
e_1 \wedge e_2, \;
e_3 \wedge e_2 + e_1 \wedge e_4 + p \, e_3 \wedge e_4
\rangle \leq W.
\]
Note that $L \otimes \QQ_p$ is generated by the vectors $e_1 \wedge e_2$ and
$(e_1 + p e_3) \wedge (e_2 + p e_4)$, and so $\SK_1(\Lie \otimes \QQ_p) = 0$.
The same two vectors do not generate $L$. In fact, the only decomposable
wedges in $L$ are $\ZZ_p$-multiples of $e_1 \wedge e_2$ and $(e_1 + p e_3)
\wedge (e_2 + p e_4)$. Therefore $\SK_1(\Lie) \cong C_p$.
\end{example}

Relying on Lemma \ref{l:sk1_subspaces_submodules}, it is therefore possible to completely
linearise the question of whether or not $\SK_1(\ZZ_p \llbracket G \rrbracket)$
is finite, i.e., we are reduced to analysing whether or not $\SK_1(\grasslie(L)
\otimes \QQ_p)$ is trivial for a generic subspace $L \subseteq V \wedge V$,
and this is in turn a question about decomposability in the exterior algebra
$V \wedge V$ over $\QQ_p$. We
emphasize that we are dealing with nilpotent algebras. This complements
the existing literature where it is shown that kernels of commutator maps for
semisimple Lie algebras are always generated by decomposable elements, both over
$\QQ_p$ (see \cite[Theorem 3.1]{SchVen13lie}) and $\CC$ (see  \cite[Corollary
5.1]{Kos65}).

\section{The decomposability map over $\QQ_p$}
\label{section:decomposability_Qp}


\noindent
Our aim in this section is to analyse the behaviour of $\SK_1$ in the generic case.
Based on the previous chapter, this boils down to understanding the points
in the Grassmannian $\Gr(r, W)$ that can be generated by decomposable elements
$D \Lie \cap \ker_{\Lie}$. We therefore establish a decomposability map 
whose image are precisely the points for which the Bogomolov multiplier
of the objects associated to it under our parametrization vanishes.
We study the local properties of this map and give a precise
description of its differential. The bulk of the work is devoted to finding
points in which the differential is a local isomorphism under a certain
restriction on the parameters of the construction.

\subsection{The decomposability map} \label{ss:decomposability_morphism_Qp}

Let $V$ be a vector space over $\QQ_p$ of dimension $d$.
Fix some $1 \leq r \leq d$. The variety $\Gr(r, W)$ consists of vector
subspaces of $W = V \wedge V$ of dimension $r$. The subspaces that can be 
generated by decomposable wedges have a basis consisting of $r$
decomposable wedges. The set of decomposable wedges in $W$ itself forms
a variety $\Gr(2, V)$. 
We therefore have a rational map of projective varieties
\begin{equation} \label{e:definition_of_Psi_Qp}
\textstyle
\Psi \colon \Gr(2, V)^r \dashrightarrow \Gr(r, W), \quad 
(L_1, \dots, L_r) \mapsto 
\bigwedge^2 L_1 \oplus \dots \oplus \bigwedge^2 L_r.
\end{equation}
Note that $\Psi$ is not regular everywhere, since the spaces $L_i$ might
overlap. It is, however, defined on a Zariski dense open subset. 
The subspaces of $W$ that can be
generated by decomposable wedges, i.e., those whose corresponding Lie algebra
has a trivial $\SK_1$, are precisely the elements of $\im \Psi$.

The map \eqref{e:definition_of_Psi_Qp} can be expressed in
terms of Pl\"ucker coordinates on the Grassmannians (see \cite{Sch94}),
giving a more explicit map
\[
\Psi \colon \Gr(2, V)^r \dashrightarrow \Gr(r, W), \quad 
(x_1 \wedge y_1, \dots, x_r \wedge y_r) \mapsto 
\langle x_1 \wedge y_1, \dots, x_r \wedge y_r \rangle
\]
with
\[
\Gr(2, V) \subseteq \PP^{\binom{d}{2} - 1}, \qquad
\Gr(r, W) \subseteq \PP^{\binom{\binom{d}{2}}{r} - 1}.
\]
We call $\Psi$ the {\em decomposability map}.
Observe that the dimension of the domain of $\Psi$ is equal to
$r \cdot 2(d-2)$, while the dimension of the codomain of $\Psi$
is equal to $r ( \binom{d}{2} - r )$.
The difference between the two is
\[
\dim (\codom \Psi) - \dim (\dom \Psi) =  r \left( \binom{d-2}{2} + 1 - r \right).
\]

\subsection{Few relators}

In order to understand the map $\Psi$ better, we first deal with the case when
the dimension $r$ of the selected subspace is not too close to the full
dimension $\binom{d}{2}$.

\begin{lemma} \label{l:grassmanian_small_r}
Suppose $r \leq \binom{d-2}{2}$. Then $\im \Psi$ is contained in a proper
subvariety of $\Gr(r, W)$.
\end{lemma}
\proof
The inequality in the statement is equivalent to saying that the dimension of
the domain of $\Psi$ is strictly smaller than the dimension of its codomain.
\endproof

Therefore a generic subspace of a small dimension $r$ produces a Lie algebra
whose $\SK_1$ is not trivial, and therefore a group whose $\SK_1$ is not finite.
In other words, fixing the dimension $r$ and letting the number of generators $d$
tend to infinity, a generic group will have an infinite $\SK_1$.

\subsection{Plenty of non-decomposable subspaces} \label{ss:plenty_of_nondecomposable_subspaces}

In general, the image $\im \Psi$ will not be contained in a proper subvariety of
$\Gr(r, W)$. However, we can (almost) always find an asymptotically large
subset of $\Gr(r, W)$ that does not intersect $\im \Psi$.

\begin{proposition} \label{p:large_subvariety_Qp}
Suppose that $\rho = \binom{d}{2} - r \geq 5$.
Then there is a subvariety $Q_r \subseteq \Gr(r, W)$ such that
\[
Q_r \cap \im \Psi= \emptyset
\qquad \text{and} \qquad
\lim_{\rho \to \infty}\frac{\dim Q_r}{\dim \Gr(r,W)} = 1.
\]
\end{proposition}
\proof

Let $\{ v_i \mid 1 \leq i \leq d \}$ be a basis of $V$ and let $\bar V 
= \langle v_1, v_2, v_3, v_4 \rangle$. Set $X = \langle x \rangle
 \leq \bar V \wedge \bar
V$ to be a $1$-dimensional subspace that is not generated by
decomposable wedges, for example $x = v_1 \wedge v_2 + v_3 \wedge v_4$.
Now let $U \leq W$ be the standard complement of
$\bar V \wedge \bar V$, so that  $v_k \wedge v_l \in U$ as
long as not both $v_k, v_l$ belong to $\bar V$. Thus 
$W = (\bar V \wedge \bar V) \oplus U$.

Fix a subspace $L = \langle \ell_1, \dots, \ell_{r-1} \rangle \leq U$ of
dimension $r - 1$ and its complement $K$, so that $U = L \oplus K$. Now select
a set of $r$ vectors $k_0, \dots, k_{r-1}$ in $K$ and associate to it the vector
subspace
\[
S = \langle x + k_0, \ell_1 + k_1, \dots, \ell_{r-1} + k_{r - 1} \rangle \in \Gr(r, W).
\]
Note that $S + U = X \oplus U$, which is not generated by decomposable wedges.
Whence $S$ has the same property. 

It follows that every choice of vectors as above produces a  subspace of $W$
of codimension $\rho$ that is not generated by decomposable wedges. Moreover,
these subspaces form a quasiprojective subvariety of the Grassmannian, call it
$Q_r$. This subvariety is contained in the open subset given by the non-vanishing
of the $r \times r$ minor associated to $(x, \ell_1, \dots, \ell_r)$. Inside
this open subset, $Q_r$ is determined by the vanishing of the minors associated
to selecting a basis vector of a complement of $X$ in $\tilde V \wedge \tilde
V$ and $r-1$ vectors in $(x, \ell_1, \dots, \ell_r)$.

Counting dimensions, we have $\dim Q_r = r (\binom{d}{2} - 6 - (r-1))$ and so
\[
\frac{\dim Q_r}{\dim \Gr(r,W)} = 1 - \frac{5}{\rho},
\]
whence the proposition.
\endproof

\subsection{The differential of decomposability} \label{ss:differential_of_decomposability}

We now proceed with our analysis of $\Psi$. Let $\mathbf{L} = (L_1, \dots,
L_r)$ be a point in the domain of $\Psi$. We are interested in the behaviour
of $\Psi$ around $\mathbf{L}$, so we will determine its differential
\[
\der \Psi_{\mathbf{L}}
\colon
T_{\mathbf{L}} \Gr(2,V)^r
\to
T_{\Psi(\mathbf{L})} \Gr(r, W).
\]

\subsubsection*{Tangent spaces of Grassmannians}

Tangent vectors $T_{L_i}\Gr(2,V)$ can be identified with linear morphisms
$\hom(L_i, V / L_i)$ in the following way (see \cite{Voi02}).
Let $A \in \hom(L_i, V)$ be a
linear map. Consider its associated curve
\[
\gamma \colon \QQ_p \to \Gr(2,V),
\quad
t \mapsto \im (id_{L_i} + t A)
\]
In coordinates, this can be expressed as follows.
Suppose $L_i = \langle x_i, y_i \rangle$.
Then
\[
\im (id_{L_i} + t A) = \langle x_i + t A x_i, y_i + t A y_i \rangle.
\]
Correspondingly, we obtain the tangent vector 
$[\gamma] \in T_{L_i} \Gr(2,V)$. This is the vector associated to the
linear map $A$. Replacing $A$ by a map $B \in \hom(L_i, V)$ gives the same
tangent vector if and only if $\im(A - B) \leq L_i$. This explains the
identification
\[
T_{L_i}\Gr(2,V) \equiv \hom(L_i, V / L_i),
\]
and the same reasoning gives the identification on the right hand side,
\[
T_{\Psi(\mathbf{L})} \Gr(r,W) \equiv \hom(\Psi(\mathbf{L}), W / \Psi(\mathbf{L})).
\]
Both of these will henceforth be used.

\subsubsection*{The matrix of the differential}

The differential $\der \Psi_{\mathbf{L}}$ is completely determined by its
values on each of the components $T_{L_i} \Gr(2,V)$ of its domain. Let us
therefore fix a linear map $A \in \hom(L_i, V)$.
Its associated curve $\gamma$ induces a curve in $\Gr(2,V)^r$ by
means of
\[
\tilde \gamma \colon t \mapsto 
(L_1, \dots, L_{i-1}, \gamma(t), L_{i+1}, \dots, L_r) \in \Gr(2,V)^r.
\]
In order to determine $\der \Psi_{\mathbf{L}}[\tilde \gamma]$,
we first map the curve $\tilde \gamma$ with $\Psi$ into
\[
\textstyle
\Psi \circ \tilde \gamma \colon \QQ_p \to \Gr(r, W),
\quad
t \mapsto
\bigwedge^2 \gamma(t)
\oplus
\bigoplus_{j \neq i} \bigwedge^2 L_j.
\]
The volume form $\gamma(t) \wedge \gamma(t)$ can be expressed in terms
of coordinates as
\[
(x_i + t A x_i) \wedge (y_i + t A y_i) =
x_i \wedge y_i +
t \left( A x_i \wedge y_i + x_i \wedge A y_i \right) +
t^2 A x_i \wedge A y_i.
\]
The tangent vector
$[\Psi \circ \tilde \gamma]$
can therefore be represented by the linear map
$\mathbf{A} \in \hom(\Psi(\mathbf{L}), W / \Psi(\mathbf{L}))$
given as
\begin{align*}
\mathbf{A} \lvert_{L_i \wedge L_i} &= 
(x_i \wedge y_i \mapsto A x_i \wedge y_i + x_i \wedge A y_i), \\
\mathbf{A} \lvert_{L_j \wedge L_j} &= 
0 \text{ for $j \neq i$}.
\end{align*}
In this notation, we therefore have $\der \Psi_{\mathbf{L}} (A) = \mathbf{A}$.

\subsection{Submersiveness of decomposability} \label{ss:submersiveness}

Recall from Lemma \ref{l:grassmanian_small_r} that as long as $r \leq \binom{d-2}{2}$, 
the image $\im \Psi$ is contained in a proper subvariety. We will
now show that for all other values of $r$, $\Psi$ is a local surjection.

\begin{theorem} \label{t:psi_submersion}
The map $\Psi$ is a submersion on a dense open subset
of its domain
if and only if $r \geq \binom{d-2}{2} + 1$.
\end{theorem}

The bound value is particularly interesting. It will follow from our proof that
when $r = \binom{d-2}{2} + 1$, the decomposability map $\Psi$ is a local
isomorphism.

First of all, it is clear that  the points in which $\Psi$ is a submersion
form a Zariski open subset of the domain of $\Psi$, since they correspond to
points in which the differential is of full rank (see \cite[Section
III.10]{Ser09}). In order to show that the principal minors of the differential are not
all identically zero, we will provide explicit points in which the
differential is indeed of full rank, and so the set of points in which $\Psi$
is a submersion is indeed non-empty and dense. This will be done by a
process of understanding more precisely the condition for $\Psi$ to
be a submersion.

\begin{remark}
It will follow from the proof that the same result holds in any 
complete field of characteristic zero, in particular $\mathbb{C}$.
\end{remark}

\subsubsection*{Affine charting the Grassmannian}

We will use the following chart on the Grassmannian $\Gr(2, V)$.
Suppose $K \leq V$ is a subspace of codimension $2$.
Choose a complement $L \leq V$ so that $V = K \oplus L$.
Any other complement can be obtained as the image of the map
\[
\varphi_K \colon \hom(L, K) \to \Gr(2,V),
\quad
f \mapsto \im(id_{L} + f).
\]
This map is an immersion of the affine space $\hom(L, K) \cong \QQ_p^{2(d-2)}$
into $\Gr(2,V)$ with a Zariski open image. We will denote the image of this
open chart by $U_K$. Note that
\[
U_K = \{ L' \in \Gr(2,V) \mid L' \cap K = 0 \},
\]
which is independent of the choice of the initial complement $L$.

\subsubsection*{Criterion for submersiveness}

We now develop a criterion on when $\Psi$ is a submersion at a
point. This relies on knowing the differential of $\Psi$. Based on Subsection
\ref{ss:differential_of_decomposability}, we have
\[
\im \der \Psi
=
\left\{
\mathbf{A} \in \hom(\Psi(\mathbf{L}), W / \Psi(\mathbf{L})) 
\; \mid \;
\forall i. \;
\im \mathbf{A} \lvert_{L_i \wedge L_i} \subseteq \Psi(\mathbf{L}) + (L_i \wedge V)
\right\}.
\]
Whence the condition for $\Psi$ to be a submersion is equivalent to 
\begin{equation} \label{eq:condition_for_submersion}
\forall i. \quad \Psi(\mathbf{L}) + (L_i \wedge V) = W.
\end{equation}
Consider the natural projections $\pi_i \colon V \to V/L_i$.
These induce maps
\[
\textstyle
\bigwedge^2 \pi_i \colon
W \to \bigwedge^2 (V/L_i),
\]
using which conditions \eqref{eq:condition_for_submersion}
can be restated as
\begin{equation} \label{eq:condition_for_submersion_with_projections}
\textstyle
\forall i. \quad 
\left(\bigwedge^2 \pi_i \right)
\left( \Psi(\mathbf{L}) \right)
=
\bigwedge^2 (V/L_i).
\end{equation}
This is the condition that we will use in what follows.
Note that since 
\[
\textstyle
\dim \Psi(\mathbf{L}) \leq r
\quad \text{and} \quad
L_i \wedge L_i 
\subseteq 
\Psi(\mathbf{L})
\cap
\ker \left( \bigwedge^2 \pi_i \right),
\]
a necessary condition for \eqref{eq:condition_for_submersion_with_projections}
is that $r - 1 \geq \binom{\dim (V/L_i)}{2} = \binom{d-2}{2}$.

\subsubsection*{Condition for submersiveness in a chart}

We will now develop the condition
$\eqref{eq:condition_for_submersion_with_projections}$ further by restricting
the subspaces $L_i$ to belong to an open subset $U_K \subseteq \Gr(2,V)$ for
some fixed subspace $K \leq V$ of codimension $2$.  Select a complement $L$ of
$K$. Each of the spaces $L_i$ therefore arises uniquely from a linear map $f_i
\colon L \to K$ via the affine chart $\varphi_K \colon \hom(L, K) \to
\Gr(2,V)$. We will take $f_1$ to be the zero map, so that $L_1$ is just the
distinguished complement $L$.

The advantage of restricting to $U_K$ is that we can replace
the quotient $V/L_i$ by the concrete model $K$ 
in each of the natural projections $\pi_i \colon V \to V/L_i$.
This corresponds to considering the alternative projections
\[
\tilde \pi_i = f_i - id_K \colon L \oplus K = V \to K
\]
with kernels $\im(id_L + f_i) = L_i$. In these terms, we can rephrase the
condition \eqref{eq:condition_for_submersion_with_projections} as
\begin{equation} \label{eq:condition_for_submersion_with_fixed_quotient}
\textstyle
\forall i. \quad 
\left(\bigwedge^2 \tilde \pi_i \right)
\left( \Psi(\mathbf{L}) \right)
=
K \wedge K.
\end{equation}
In order to determine the left hand side, we capitalize on 
the fact that the subspaces $L_i$ are determined by the maps
$f_i$. For a subspace $L_j \in U_K$, we have
\[
\textstyle
\tilde \pi_i
\left( L_j \right)
=
\im((f_i - id_K) \circ (id_L + f_j)) 
=
\im(f_j - f_i) \subseteq K.
\]
Condition \eqref{eq:condition_for_submersion_with_fixed_quotient}
is therefore saying that the volume forms of $\im(f_j - f_i)$
should generate $K \wedge K$ for any fixed $i$,
\begin{equation} \label{eq:condition_volume_forms_generate_KwK}
\textstyle
\forall i. \quad
K \wedge K
=
\sum_{j = 1}^r
\bigwedge^2 \im(f_j - f_i).
\end{equation}

\subsubsection*{Coordinates and the canonical choice}

We will now focus on the particular case when $r$ is equal to the bound value
$\binom{d-2}{2} + 1$, and resolve the submersiveness condition in this case.
Our objective is to find examples of tuples $\mathbf{L}$ that satisfy the
condition and this can be achieved, as shown below, by a fairly standard
choice of subspaces $L_i$.

We begin by introducing coordinates in the above construction. Let 
$e_1, \ldots, e_d$ be a basis of $V$, and take
\[
L = \langle e_1, e_2 \rangle,
\quad
K = \langle e_3, \dots, e_d \rangle.
\]
Note that $L_1 = L$. The remaining
$r-1 = \binom{d-2}{2}$ subspaces $L_i$ will be parametrized
by the set
\[
\II = \{
I \subseteq \{ 3, \dots, n \}
\mid
\lvert I \rvert = 2
\}
\]
in the following way.
For each $I = \{ a, b \} \in \II$ with $a < b$, we set $\omega_I = e_a \wedge
e_b \in \bigwedge^2 K$. The set $\{ \omega_I \mid I \in \II \}$ forms a
standard basis of the space $\bigwedge^2 K$.  The subspace of $K$ whose volume
form is represented by the wedge $\omega_I$ will be denoted by $\Lambda_I =
\langle e_a, e_b \rangle \subseteq K$.
We will look for the subspaces $L_I$, now parametrized by
the set $\II$, via their corresponding $f_I$ maps. 
In view of our understanding of condition \eqref{eq:condition_for_submersion_with_fixed_quotient}
via volume forms of images of the $f_I$ maps,
we will restrict ourselves to the case when these maps
are embeddings with canonical images, i.e.,
\[
\forall I \in \II. \quad
f_I \colon L = \langle e_1, e_2 \rangle \to \Lambda_I = \langle e_a, e_b \rangle
\]

Note that the set $\II$ can be seen as the vertex set of the Johnson graph
$J(n, 2)$. Two vertices $I, J$ in this graph are connected, written as
$I - J$, if and only if
$\lvert I \cap J \rvert = 1$. We will use this terminology in what follows.

\subsubsection*{Condition for submersiveness in the canonical choice}

Let us now determine the conditions on the maps $f_I$ so that
the projections $\tilde \pi_1$ and $\tilde \pi_I$ satisfy
\eqref{eq:condition_for_submersion_with_fixed_quotient}.

First of all, when projecting according to  $\tilde \pi_1 \colon V \to K$ with
kernel $L_1$, the volume forms $\bigwedge^2 \im(f_J - f_1) = \omega_J$ for $J
\in \II$ clearly generate $\bigwedge^2 K$. This takes care of the exceptional
case.

Now consider projecting according to $\tilde \pi_I$ for some $I \in \II$.
The volume forms of $\im(f_1 - f_I) = \im(f_I)$ and
of $\im(f_J - f_I)$ for all $J \in \II$ generate $\bigwedge^2 K$
if and only if
\begin{equation} \label{eq:condition_for_submersion_canonical_choice}
\textstyle
\forall H \in \II. \quad
\omega_H \in
\left\langle
\omega_I
\right\rangle
+
\sum_{J \in \II}
\bigwedge^2 \im(f_J - f_I)
\end{equation}
For every $J \in \II$, we can write the volume form of
$\im(f_J - f_I)$
expicitly as
\[
\left(
f_I(e_1 \wedge e_2) + f_J(e_1 \wedge e_2)
\right)
-
\left(
f_I(e_1) \wedge f_J(e_2) + f_J(e_1) \wedge f_I(e_2)
\right)
\]
Observe that
\[
f_I(e_1 \wedge e_2)
-
\left(
f_I(e_1) \wedge f_J(e_2) + f_J(e_1) \wedge f_I(e_2)
\right)
\in
\langle \omega_I \rangle
+
\left\langle
\left\{
\omega_H
\mid
H \in \II, H - I
\right\}
\right\rangle
\]
and that $f_J(e_1 \wedge e_2)$ is a non-trivial multiple of $\omega_J$.
This means that as long as \eqref{eq:condition_for_submersion_canonical_choice}
is shown to be satisfied for all subsets $H \in \II$ with 
that are connected to $I$, 
it will also be satisfied for all the other subsets $H \in \II$.
We can therefore restrict ourselves to inspecting the situation
when $I = \{ a, b \}$ and $H$ is either $\{ a, c \}$ or $\{ b, c \}$.
In such a case, consider the corresponding part of the space
in \eqref{eq:condition_for_submersion_canonical_choice},
\begin{equation} \label{eq:condition_on_intersecting_pairs}
\textstyle
\langle \omega_{\{ a, b \}} \rangle
+
\bigwedge^2 \im(f_{\{a,c\}} - f_{\{a,b\}})
+
\bigwedge^2 \im(f_{\{b,c\}} - f_{\{a,b\}}).
\end{equation}
Observe that the space \eqref{eq:condition_on_intersecting_pairs}
is a subspace of
$\langle \omega_{\{a,b\}}, \omega_{\{a,c\}}, \omega_{\{b,c\}} \rangle$.
These two spaces are equal if and only if
both $\{a,c\}$ and $\{b,c\}$
satisfy $\eqref{eq:condition_for_submersion_canonical_choice}$.
Whence \eqref{eq:condition_for_submersion_with_fixed_quotient}
is satisfied if and only if the vectors
\begin{equation} \label{eq:condition_independence_of_3_vectors}
\textstyle
\left\{
\omega_{I},
\;
\bigwedge^2 \im(f_{J} - f_{I}),
\;
\bigwedge^2 \im(f_{H} - f_{I})
\right\}
\end{equation}
are linearly independent for all triangles $I - J - H - I$
in the graph.

\subsubsection*{Transition maps}

In order to express this last condition more clearly,
set $f_I^{[-1]} \colon \Lambda_I \to L$ to be the inverse of $f_I$
defined on its image $\Lambda_I$, and let
\[
\delta_{J}^I \colon \Lambda_J \to \Lambda_I,
\quad
\delta_J^I = f_I \circ f_J^{[-1]}.
\]
These maps should be thought of as transition maps as will be clear
from what follows.
For all $a,b$, we can pull back the basis of $\Lambda_{\{a,b\}}$ into
\[
u_{\{a,b\}, a} = f_{\{a,b\}}^{-1}(e_a),
\quad
u_{\{a,b\}, b} = f_{\{a,b\}}^{-1}(e_b).
\]
The advantage of this new basis is that we can replace the volume forms
in \eqref{eq:condition_independence_of_3_vectors} with
their scalar multiples
\[
\textstyle
(f_{\{a,c\}} - f_{\{a,b\}})(u_{\{a,c\}, a} \wedge u_{\{a,c\},c}),
\;
(f_{\{b,c\}} - f_{\{a,b\}})(u_{\{b,c\}, b} \wedge u_{\{b,c\},c}).
\]
The first of these vectors can be rewritten as
\[
\textstyle
\left(e_a - \delta_{\{a,c\}}^{\{a,b\}}(e_a)\right)
\wedge
\left(e_c - \delta_{\{a,c\}}^{\{a,b\}}(e_c)\right),
\] 
which is simply
\[
\textstyle
\left( e_a - \delta_{\{a,c\}}^{\{a,b\}}(e_a) \right) \wedge e_c
\pmod{\langle \omega_{\{a,b\}} \rangle}.
\]
We have a similar expression for the second vector. The vectors \eqref{eq:condition_independence_of_3_vectors} are therefore independent if and only if
the vectors
\begin{equation} \label{eq:condition_delta_2_vectors}
\left\{
e_a - \delta_{\{a,c\}}^{\{a,b\}}(e_a),
e_b - \delta_{\{b,c\}}^{\{a,b\}}(e_b)	
\right\}
\end{equation}
are independent.

Let $X = \{ 3, 4 \} \in \II$. In order to resolve the last condition, 
Note that we can express $\delta_J^I = (\delta_I^X)^{-1} \circ \delta_J^X$.
The subfamily
$\Delta^X = \left\{ \delta_{I}^X \mid I \in \mathcal{I} \right\}$
of maps into $\Lambda_X$ therefore uniquely
determines all the $\delta_J^I$ maps.
Moreover, we can apply $\delta_{\{a,b\}}^X$ to convert
the vectors \eqref{eq:condition_delta_2_vectors} into
\begin{equation} \label{eq:condition_delta_2_symmetric_vectors}
\left\{
\left( \delta_{\{a,b\}}^X - \delta_{\{a,c\}}^X \right)(e_a),
\left( \delta_{\{a,b\}}^X - \delta_{\{b,c\}}^X \right)(e_b),
\right\}
\end{equation}
obtaining a condition expressed solely in terms of the maps
belonging to $\Delta^X$.

It is possible to invert the described procedure and recover the maps $f_I$
from the transition maps $\delta_J^X$ and the map $f_X$. To achieve this, we
can simply define $f_I = \left(\delta_I^X\right)^{-1} \circ f_X$ for any $I
\in \II$. We can even assume $f_X$ is given as the standard map $e_1 \mapsto
e_3$, $e_2 \mapsto e_4$. Our problem is thus reduced to finding a family of
maps $\Delta^X$ such that the vectors
\eqref{eq:condition_delta_2_symmetric_vectors} are independent.

\subsubsection*{Finding the family $\Delta^X$}

We can take $\delta_X^X = id_{\Lambda_X}$
We search for the remaining $\delta_I^X$ maps
in terms of their matrices written in
standard basis of $\Lambda_I$ and $\Lambda_X$. Set
\[
\delta_I^X \equiv
\begin{pmatrix}
\alpha_I & \gamma_I \\ 
\beta_I & \delta_I
\end{pmatrix} \in M_2(\QQ_p).
\]
The maps $\delta_I^X$ have to be invertible, giving the conditions
$\alpha_I \delta_I - \beta_I \gamma_I \neq 0$
for all $I \in \mathcal{I}$.
The condition that the vectors \eqref{eq:condition_delta_2_symmetric_vectors}
be independent can be written as
\[
\forall (I,J,H) \in \II^3. \quad
I - J - H - I
\;
\Longrightarrow
\;
\det \begin{pmatrix}
\alpha_I - \alpha_J & \gamma_I - \gamma_H \\
\beta_I - \beta_J & \delta_I - \delta_H
\end{pmatrix} \neq 0.
\]
The set of solutions to these conditions forms a complement of
a finite union of quadrics in $M_2(\QQ_p)$, and so the set of
solutions to our problem is Zariski dense. To write down an
explicit solution, we can simply take scalar matrices
with the conditions that $\alpha_I \neq 0$ and $\alpha_I \neq \alpha_J$
for $I - J$. 
Taking this back to the subspaces $L_i$, we obtain the following 
quite elementary example.

\begin{example} \label{ex:point_of_local_isomorphism}
Let $V = \langle e_1, \dots, e_n \rangle$.
Set $L_1 = \langle e_1, e_2 \rangle$ and
for any $3 \leq i < j \leq n$,
set
$L_{i,j} = 
\left\langle 
e_1 + \lambda_{i,j} e_i, \;
e_2 + \lambda_{i,j} e_j 
\right\rangle$
for non-zero distinct scalars $\lambda_{i,j} \in \QQ_p$.
Then
$\Psi$ is a submersion in a neighbourhood of $\mathbf{L}$.
\end{example}

\subsubsection*{Larger than bound}

In order to complete the proof of Theorem \ref{t:psi_submersion},
it remains to inspect the case when $r > \binom{d-2}{2} + 1$. 
To this end, let $Z$ be any subset of $\{ 1, \dots, r \}$ with
$|Z| = \binom{d-2}{2} + 1$. Consider the coordinate projection
\[
\pr_Z \colon \Gr(2,V)^r \to \Gr(2,V)^{\binom{d-2}{2} + 1}
\]
onto the $Z$-axes. We have shown above that there is a
dense neighbourhood
$\mathcal{U}_Z \subseteq \Gr(2,V)^{\binom{d-2}{2} + 1}$
consisting of those tuples $\mathbf{L}$ in which the associated
map $\Psi$ is a submersion, i.e., satisfies the condition
\eqref{eq:condition_volume_forms_generate_KwK}.
Let
\[
\mathcal{U} = 
\bigcap_{
\begin{smallmatrix}
Z \subseteq \{ 1, \dots, r \} \\
|Z| = \binom{d-2}{2} + 1
\end{smallmatrix}
}
\pr_Z^{-1}(\mathcal{U}_Z)
\subseteq
\Gr(2,V)^r.
\]
Note that $\mathcal{U}$ is a Zariski dense open subset of $\Gr(2,V)^r$.
Clearly every point $\mathbf{L} \in \mathcal{U}$ satisfies \eqref{eq:condition_volume_forms_generate_KwK}. The proof of Theorem
\ref{t:psi_submersion} is thus complete.

\subsection{Immersiveness of decomposability} \label{ss:immersiveness}

We now deal with complementing the previous section by proving the following.

\begin{theorem} \label{t:psi_immersion}
The map $\Psi$ is an immersion on a dense open subset
of its domain
if and only if $r \leq \binom{d-2}{2} + 1$.
\end{theorem}

Our proof will rely on exploiting that in the bound case when
$r = \binom{d-2}{2} + 1$, we already know that $\Psi$ is
a local isomorphism on a dense open subset. Recall that as before, it suffices
to find one good point of the domain of $\Psi$. To this end, we connect
the general case with the bound one. First of all,
consider the coordinate projection
\[
\pr \colon \Gr(2, V)^{\binom{d-2}{2} + 1} \to \Gr(2,V)^r
\]
onto the first $r$-axes. Both the source and the target of
$\pr$ have an associated rational map $\Psi$,
\[
\begin{tikzcd}
\Gr(2,V)^{\binom{d-2}{2} + 1} \arrow[r, "\pr"] \arrow[d, dashed, "\Psi"] &
\Gr(2,V)^r \arrow[d, dashed, "\Psi"] \\
\Gr(\binom{d-2}{2} + 1, W) &
\Gr(r, W).
\end{tikzcd}
\]
These induce differential maps,
\[
\begin{tikzcd}
T_{\mathbf{L}}\Gr(2,V)^{\binom{d-2}{2} + 1} 
\arrow[r, two heads, "\der\pr_{\mathbf{L}}"] 
\arrow[d, "\der\Psi_{\mathbf{L}}"] &
T_{\pr(\mathbf{L})}\Gr(2,V)^r 
\arrow[d, "\der\Psi_{\pr(\mathbf{L})}"] 
\arrow[l,out=150,in=30,dotted,"\sigma"]
\\
T_{\Psi(\mathbf{L})}\Gr(\binom{d-2}{2} + 1, W) \arrow[r, dotted,"\iota"] &
T_{\Psi(\pr(\mathbf{L}))}\Gr(r, W).
\end{tikzcd}
\]
Note that $\der\pr_{\mathbf{L}}$ is surjective
and has a natural splitting $\sigma$.
As long as $\mathbf{L}$ belongs to some dense open subset
$\mathcal{U} \subseteq \Gr(2,V)^{\binom{d-2}{2} + 1}$,
the map $\der\Psi_{\mathbf{L}}$ is an isomorphism.
Select complements $W = \Psi(\mathbf{L}) \oplus X$
and $\Psi(\mathbf{L}) = \Psi(\pr(\mathbf{L})) \oplus Y$.
After identifying the tangent spaces with 
hom sets, the map $\iota$ is the natural composition
\[
\iota \colon 
\hom(\Psi(\mathbf{L}), X)
\xrightarrow{res}
\hom(\Psi(\pr(\mathbf{L})), X)
\to
\hom(\Psi(\pr(\mathbf{L})), X \oplus Y)
\]
Now, as $\iota$ is an embedding on the image of 
$\der\Psi_{\mathbf{L}} \circ \sigma$, it follows that 
$\der\Psi_{\pr(\mathbf{L})}$ is injective.
This completes the proof of Theorem \ref{t:psi_immersion}.
Example \ref{ex:point_of_local_isomorphism}
gives an explicit point $\mathbf{L} \in \mathcal{U}$, from which
points $\pr(\mathbf{L})$ can be produced.

\begin{remark}
One can think of searching for a point $\mathbf{L}$
in which $\Psi$ is a local immersion
by searching for points over which the map
$\Psi$ has a $0$-dimensional fibre.
This means that we can pass to the target of $\Psi$
and search there.
Consider $\Psi(\mathbf{L})$ as a 
subspace of $W$ of dimension $r$.
This subspace is spanned by
its intersections with
the image of the Plücker embedding of $\Gr(2,V)$ into $W$,
\[
\Psi(\mathbf{L}) =
\langle
\Psi(\mathbf{L}) \cap \Gr(2, V)
\rangle.
\]
The latter condition in fact characterizes $\im \Psi$,
since for $H \in \Gr(r,W)$ we have
$H \in \im \Psi$ if and only if $H = \langle H \cap \Gr(2,V) \rangle$.
Our search for the $r$ points therefore corresponds
to finding a linear section of the Grassmannian $\Gr(2,V)$
that is $0$-dimensional and is generated
by its finitely many points. 
Over the algebraic closure $\bar \QQ_p$,
one can use Bertini's theorem
(see \cite[Theorem 12.1 (i) and (ii)]{Arr17}) together with B\'ezout to deduce
that generic subspaces have this property.
More precisely, there is a dense open subset
$\mathcal{V} \subseteq \Gr(\binom{d-2}{2} + 1, d)$
such that every element $\Lambda \in \mathcal{V}$
has the property that $\dim ( \Lambda \cap \Gr(2,V)) = 0$.
In this language, our Example \ref{ex:point_of_local_isomorphism}
shows that such properties can be obtained from
rational points both in
$\Gr(\binom{d-2}{2} + 1, d)$
and
$\Gr(2,V)$.
\end{remark}



\section{The decomposability map over $\FF_p$}
\label{section:decomposability_Fp}


We now descend from groups and algebras over $\ZZ_p$ to their finite
quotients. The role of torsion-freeness will be played here by the assumption
that the finite groups are of exponent $p$. All the structure of these groups
is therefore captured in the relations between commutators. We apply the same
approach as above to set up a parametrization of these groups, express their
Bogomolov multipliers and define the decomposability map. We first
illustrate what happens with a small number of generators, and inspect when
this map is surjective. We analyse the generic behaviour for large primes
$p$  by using the results in characteristic $0$.  On the other hand, when
fixing $p$, we show how the subvariety of large dimension representing groups
with non-trivial Bogomolov multipliers from the previous section can be made to
exist in the finite case as well. We conclude with an application of our
technique regarding a more elementary problem in commutators in finite groups,
which can also be expressed in the language of Grassmannians. Throughout this
section, we will prefer to use the more standard term Bogomolov multiplier
rather than $\SK_1$.

\subsection{Groups and the Grassmannian variety} \label{ss:groups_and_grassmannian}

Let $G$ be a finite $p$-group of exponent $p$ (assume $p > 2$) and nilpotency
class $2$. Suppose $G$ is of rank $d$, so that $V = G / [G,G] \cong \FF_p^d$. The
structure of $G$ is then completely determined by the set of relations between
its commutators. These form a certain linear subspace in $V \wedge V$. Fix the
dimension $r$ of this space of relations. Thus we are looking at
$r$-dimensional subspaces of the $\binom{d}{2}$ dimensional $\FF_p$-space $W =
V \wedge V$. These form the Grassmannian variety $\Gr(r,W)$, and each of its
points $L$ determines a group $G_L$ obtained by imposing precisely the
relations of the corresponding subspace $L \leq W$.
This is a mod $p$ version of the parametrization \eqref{eq:parametrization_Qp},
so we stay with the same notation as in the previous parts.

\subsection{Bogomolov multiplier}

Let $L \in \Gr(r,W)$ be a subspace. The Bogomolov multiplier  of the
corresponding group $G_L$ can be recognized as follows (see \cite{Bog87}).
Let $D W$ denote the set of decomposable wedges of $W$.
These form a variety $\Gr(2,V)$.
Then we have
\[
\B_0(G_L) = \frac{L}{\langle D W \cap L \rangle}.
\]
Therefore deciding whether or not the Bogomolov multiplier is trivial
reduces to deciding whether or not the subspace $L$ is generated by decomposable
wedges. This is just the same as in the previous section, the only difference
being that we are now considering vector spaces over $\FF_p$.

\subsection{Decomposable subspaces} \label{ss:decomposability_morphism_Fp}

In order to study the subspaces generated by decomposable wedges, we consider,
as in \eqref{e:definition_of_Psi_Qp}, the rational map
\[
\psi \colon \Gr(2,V)^r \dashrightarrow \Gr(r,W)
\]
mapping an $r$-tuple of decomposable vectors of $W$ into their span in $W$. The
locus of indeterminacy consists of the tuples which do not span an
$r$-dimensional subspace of $W$. The image of $\psi$ consists precisely of the
subspaces of $W$ whose corresponding groups have trivial Bogomolov
multipliers.

\subsection{Bounds for sizes of Grassmannians}

The size of $\Gr(k, n)$ is the
number of subspaces of dimension $k$ in an $n$-dimensional vector spaces
over a finite field with $p$ elements. This is equal to the $p$-binomial coefficient
\[
\binom{n}{k}_p = \frac
{(p^n - 1)(p^n - p) \cdots (p^n - p^{k-1})}
{(p^k - 1)(p^k - p) \cdots (p^k - p^{k-1})}.
\]
We will require the following straightforward bounds for these coefficients.

\begin{lemma} \label{l:qbinomial_bounds}
Let $k \leq n$. Then
\[
p^{k(n-k) - k} =
p^{(n-1)k - k^2} \leq \binom{n}{k}_p \leq p^{nk - (k-1)k}
= p^{k(n-k) + k}.
\]
\end{lemma}
\proof
In the definition of the binomial coefficient, use the bounds
$p^{n-1} \leq p^n - p^i \leq p^n$ and collect.
\endproof

We will also need the following bound. It is stronger than the one
above, but it includes implied constants and will only be used for
asymptotics when $p$ tends to infinity.

\begin{lemma} \label{l:grassmannian_precise_bound}
Let $k \leq n$. Then
\[
|\Gr(k,n)| = p^{k(n-k)} + O(p^{k(n-k) - 1}),
\]
where the implied constant is independent of $p$.
\end{lemma}
\proof
Immediate from \cite[Section 13.5, Theorem 6]{BF72}.
\endproof

\subsection{Few generators} \label{ss:few_generators}

Here we inspect the behaviour of the map $\psi$ when the dimension of $V$
is as small as possible. First we show that $\psi$ is surjective for $d \leq 3$.

\begin{proposition}
Suppose $d \leq 3$. Then $\psi$ is surjective.
\end{proposition}
\proof 
We are claiming that $\B_0(G_L) = 0$ for all $L \leq W$.
This is clear for $d \leq 2$, since $V \wedge V$ is at most $1$-dimensional.
When $d = 3$, we either have that $L = 0$, and so $\B_0(G_L) = 0$,
or $L \neq 0$, in which case $|G_L| \leq p^5$. As $G_L$ is also of
nilpotency class $2$, it follows from \cite{Mor12p5} that $\B_0(G_L) = 0$.
\endproof

The smallest interesting case is when $d = 4$. We analyse it in detail
by exploiting the presence of the action of $\GL(V)$.

\begin{proposition} \label{p:d=4}
Suppose $d = 4$.
\begin{enumerate}
	\item Let $r = 1$. 
	Then $\psi$ is injective and
	\[
	\lim_{p \to \infty} \left(
	\frac{|\im \psi|}{|\Gr(1,W)|} \left/  \frac1p  \right. \right)
	= 1.
	\]

	\item Let $r = 2$.
	Then $\psi$ is ``generically a $2:1$ map'' and
	\[
	\lim_{p \to \infty}
	\frac{|\im \psi|}{|\Gr(2,W)|}
	= \frac12.
	\]

	\item Let $r = 3$.
	Then $\psi$ is not surjective, but has a dense image in the sense that
	\[
	\lim_{p \to \infty} \left(
	\frac{|\im \psi|}{|\Gr(3,W)|} \left/ \left( 1 - \frac1p \right)  \right. \right)
	= 1.
	\]	
	\item Let $r \geq 4$.
	Then $\psi$ is surjective.
\end{enumerate}
\end{proposition}
\proof
{\bf (1)} Directly compute
\[
\frac{|\im \psi|}{|\Gr(1,W)|} =
\frac{\binom{4}{2}_p}{\binom{6}{1}_p} =
\frac{p^2+1}{p^3+1}.
\]

\smallskip

{\bf (2)} Set $X = \{ (L_1, L_2) \in \Gr(2,V)^2 \mid L_1 \cap L_2 = 0 \}$.

We first claim that $\psi$ restricted to $X$ has fibres of size $2$. Indeed,
let $(L_1, L_2) \in X$ with $L_1 = \langle v_1, v_2 \rangle$ and $L_2 =
\langle v_3, v_4 \rangle$. Then $\{ v_1, v_2, v_3, v_4 \}$ is a basis of $V$
and the only decomposable vectors of the space $\psi(L_1,L_2) = \langle v_1
\wedge v_2, v_3 \wedge v_4 \rangle \leq W$ are the scalar multiples of $v_1
\wedge v_2$ and  $v_3 \wedge v_4$. Each one of these lines in $W$ uniquely
determines a subspace of dimension $2$ in $V$. The only elements in $X$ that
map onto $\langle v_1 \wedge v_2, v_3 \wedge v_4 \rangle$ are therefore $(L_1,
L_2)$ and $(L_2, L_1)$.

Let us now compute the size of the set $X$. The group $\GL(V)$ acts on $W$ and
therefore on the set of its subspaces $\Gr(2,W)$. This action can be
restricted to the set $X$, on which $\GL(V)$ acts transitively. The stabilizer
of a point, say $(L_1, L_2)$, is equal to $\GL(L_1) \times \GL(L_2)$. 
In matrix terms, this corresponds to a sum of two invertible block matrices.
It follows that
\[
|X| = \left|\GL(V)\right| / | \GL\left(\FF_p^2\right) \times \GL\left(\FF_p^2\right)|
= p^4 \left(p^2+1\right) \left(p^2+p+1\right).
\]

On the other hand, we can directly compute the size of the domain of $\psi$; it is equal to
\[
\left|\Gr(2,V)^2\right| = \binom{d}{2}_p^2 = \left(p^2+1\right)^2 \left(p^2+p+1\right)^2.
\]
Note that $|X|$ is asymptotically comparable to $\left|\Gr(2,V)^2\right|$,
so that $X$ is in fact the orbit of a generic point.

We now consider the image of $\psi$. Since the restriction of $\psi$ to $X$ has fibres
of size $2$, we have 
\begin{equation} \label{e:size_of_im_psi_compared_to_x}
\frac12 |X| \leq |\im \psi| \leq (|\Gr(2,V)^2| - |X|) + \frac12 |X|.
\end{equation}
We can directly compute the size of the codomain of $\psi$,
\[
|\Gr(2,W)| = \binom{6}{2}_p = \left(p^2-p+1\right) \left(p^2+p+1\right) \left(p^4+p^3+p^2+p+1\right).
\]
Dividing inequality \eqref{e:size_of_im_psi_compared_to_x} by $|\Gr(2,W)|$, we get:
\[
\frac{p^6+p^4}{2 \left(p^6+p^4+p^3+p^2+1\right)}
\leq
\frac{|\im \psi|}{|\Gr(2,W)|}
\leq
\frac{p^6+2 p^5+5 p^4+4 p^3+6 p^2+2 p+2}{2 \left(p^6+p^4+p^3+p^2+1\right)}.
\]
Both sides converge to $\frac12$ as $p \to \infty$ and the claim is proved.

\smallskip

{\bf (3)}
In this case, $|G_L| = p^7$. 
Suppose that $B_0(G_L) \neq 0$.
It follows from \cite[Corollary 2.14, Theorem 2.13]{JezMor15unicp}
that there are only two possibilities up to isoclinism for such a group. When restricted to
groups of exponent $p$, these are in fact the only two possibilities up to isomorphism.
Therefore each one of these corresponds to an orbit of $\GL(V)$ acting on subspaces of
dimension $3$ in $W$. We analyse both cases in more detail.

The first group $G_1$ is determined by the subspace
\[
A = \langle v_1 \wedge v_2 - v_3 \wedge v_4, v_2 \wedge v_4 - \omega v_1 \wedge v_3, v_1 \wedge v_4 \rangle \leq W,
\]
where $\omega$ is a generator of $\FF_p^*$.
The stabilizer of $B$ under the action of $\GL(V)$ contains the set of matrices
of the form
\[
\left(
\begin{array}{cccc}
 a  & * & * & b  \omega  \\
 0 & c  & d  \omega  & 0 \\
 0 & d  & c  & 0 \\
 b  & * & * & a  \\
\end{array}
\right),
\]
where $0 \neq (a^2 - \omega b^2)(c^d - \omega d^2)$ and $*$ is any element of $\FF_p$. 
For any $b$, there are at most two roots of the equation $a^2 - \omega b^2 = 0$, so there are
at least $p-2$ options for selecting $a$. The same is true for the pair $c,d$.
The stabilizer therefore contains at least
\[
\left( p (p - 2) \right)^2 \cdot p^4
\]
elements, and so the orbit of $A \in \Gr(3,W)$ is of size at most
\[
\frac{|\GL(V)|}{(p - 2)^2 p^6} =
\frac{(p-1)^4 (p+1)^2 \left(p^2+1\right) \left(p^2+p+1\right)}{(p-2)^2}.
\]

The second group $G_2$ is determined by the subspace
\[
B = \langle v_1 \wedge v_2 - v_3 \wedge v_4, v_1 \wedge v_3, v_1 \wedge v_4 \rangle \leq W.
\]
The stabilizer of $B$ under the action of $\GL(V)$ contains the set of matrices
of the form
\[
\left(
\begin{array}{cccc}
 a & * & * & * \\
 0 & b & 0 & 0 \\
 0 & * & x & y \\
 0 & * & z & w \\
\end{array}
\right),
\]
where $0 \neq ab = xw - zy$ and $*$ is any element of $\FF_p$. 
There are $(p-1)^2$ options for selecting $a,b$, and after these there are $(p^2 - 1)p$
options for selecting $x,y,z,w$.
The stabilizer therefore contains at least
\[
(p - 1)^2 \cdot (p^2 - 1)p \cdot p^5
\]
elements, and so the orbit of $B \in \Gr(3,W)$ is of size at most
\[
\frac{|\GL(V)|}{(p - 1)^2 (p^2 - 1) p^6} =
(p-1) (p+1) \left(p^2+1\right) \left(p^2+p+1\right).
\]

Taking both orbits into account, we conclude that the number of subspaces of $W$
that are not in the image $\im \psi$ is at most the sum of the two bounds derived
above, which is proportional to $p^8$ as $p \to \infty$. 

On the other hand, the number of all subspaces $\Gr(3,W)$ is equal to
\[
|\Gr(3,W)| = \binom{6}{3}_p =
(p+1) \left(p^2+1\right) \left(p^2-p+1\right) \left(p^4+p^3+p^2+p+1\right),
\] 
which is proportional to $p^9$ as $p \to \infty$.

We therefore obtain
\[
\frac{|\im \psi|}{|\Gr(3,W)|} \geq 1 - a_p,
\]
where $a_p$ is a sequence proportional to $1/p$. The claim follows.

\smallskip

{\bf (4)}
In this case, $|G_L| \leq p^6$. It follows from \cite[Corollary
2.14]{JezMor15unicp} that there are no such $p$-groups of nilpotency class $2$
with non-trivial Bogomolov multipliers.
\endproof

\begin{remark}
It is clear from the proof that the action of $\GL(V)$ could be exploited
precisely because of the restriction on the number of generators.
As $d$ grows, the orbits become proportionately too small and the action
does not make a difference. We will see this explicitly
 in the last proofs of this section.
\end{remark}

\subsection{Surjectivity of decomposability} \label{ss:surjectivity}

We record a consequence of the proof of Proposition \ref{p:d=4} revealing how
the map $\psi$ is almost never surjective.

\begin{proposition} \label{p:surjectivity_iff}
The map $\psi$ is surjective if and only if $\binom{d}{2} - r \leq 2$.
\end{proposition}
\proof
Let $G_0$ be one of the $p$-groups from \cite[Theorem 2.13]{JezMor15unicp}.
The group $G_0$ is of order $p^7$, generated by $4$ elements and 
$\B_0(G_0) \neq 0$. Put $G = G_0 \times C_p^{d-4}$ for any
$d \geq 4$. Then $G$ is a $d$-generated group of exponent $p$ and
nilpotency class $2$ with $\B_0(G) \neq 0$. Hence not every subspace
of $W$ of codimension $3$ is generated by decomposable wedges, and so
$\psi$ is not a surjection for $r \leq \binom{d}{2} - 3$.

On the other hand, assume that $r \geq \binom{d}{2} - 2$ and let $L \leq W$ be
of codimension at most $2$. Suppose that $L$ is not generated by decomposable
elements. Then the set of its decomposable elements is contained in a subspace
$\bar L \leq L$ of codimension $1$. Suppose that $w \in L - \bar L$. Now, the
group $G_{\bar L}$ is of nilpotency class $2$ with $|G_{\bar L}'| \leq p^3$.
It follows from \cite[Theorem A]{Rod77} that every element of $G_{\bar L}'$ is
a commutator, and so in particular $w$ is decomposable modulo $\bar L$.
Therefore there is a decomposable vector in $L$ not belonging to $\bar L$.
This is a contradiction. We conclude that $L$ is generated by decomposable
elements, and so $\psi$ is indeed surjective.
\endproof

\subsection{Asymptotics of decomposability} \label{ss:asymptotics_of_decomposability}

We now turn to analysing the asymptotic behaviour of $\psi$. This is strongly correlated
with submersiveness and immersiveness of $\Psi$ over $\QQ_p$, although
not entirely equivalent (see Example \ref{ex:torsion_in_sk1}).

\subsubsection*{Few relators}

As long as the subspace $L$ is not of a very small codimension in $W$, the
generic Bogomolov multiplier is trivial. This is an easy application of
counting rational points, whose asymptotic behaviour coincides with
that over closed fields.

\begin{proposition} \label{p:small_r_asymptotics}
Suppose $r \leq \binom{d-2}{2}$.
Then proportion of elements in the image of  $\psi$ tends to $0$ as $p \to \infty$.
Setting $\delta = \binom{d-2}{2} - r - 3$,
we moreover have
\[
\frac{|\im \psi|}{|\Gr(r,W)|} \leq p^{- \delta r}.
\]
\end{proposition}
\proof
The asymptotic part follows immediately from Lemma 
\ref{l:grassmannian_precise_bound} applied to the source and target of $\psi$.
As for the second part, we can use the exact
upper and lower bounds from Lemma \ref{l:qbinomial_bounds} to obtain
\[
|\Gr(2,V)^r| \leq p^{2 r (d-1)}
\quad \textrm{and} \quad 
|\Gr(r,W)| \geq p^{\left( \binom{d}{2} - 1 \right) r - r^2}.
\]
It follows that
\[
\frac{|\im \psi|}{|\Gr(r,W)|} \leq p^{r^2 - r \frac{d^2 - 5 d + 2}{2}} = p^{-\delta r}.
\qedhere
\]
\endproof

\subsubsection*{Many decomposable subspaces}

We now focus on the case when $r \geq \binom{d-2}{2} + 1$. In such a
situation, the map $\Psi$ over $\QQ_p$ is a submersion almost everywhere,
meaning that its image is not contained in a proper subvariety. These differential
techniques can not help us directly with understanding the situation over
$\FF_p$, since small neighbourhoods of points over $\QQ_p$ project
entirely to a single point in $\FF_p$.
In order to transfer results over $\QQ_p$
to the case of finite fields, we will need to consider
both maps $\psi$ and $\Psi$ as particular instances of a scheme rational map
$\mathbf{\Psi}$. We will succeed in proving that the image of $\psi$ is very
large indeed.

\begin{theorem} \label{t:large_r_asymptotics}

Suppose $r \geq \binom{d-2}{2} + 1$. Then
\[
\liminf_{p \to \infty} \frac{|\im \psi|}{|\Gr(r,W)|} \geq
\left( \frac{1}{C_{d-2}} \right) ^r,
\]
where $C_{d-2}$ is the Catalan number.
In particular, the proportion of elements in the image of $\psi$
is bounded away from $0$ as $p \to \infty$.

\end{theorem}

\subsubsection*{Grassmannian schemes and the decomposability map}

To begin with, consider the schemes
$\Gr_{\mathbb{Z}}(2, d)^r$
and
$\Gr_{\mathbb{Z}}(r, \binom{d}{2})$
over $\Spec \mathbb{Z}$ (see \cite{EisHar00}).
The rational map \eqref{e:definition_of_Psi_Qp} 
is defined over $\mathbb{Z}$ and extends to a
rational map of schemes
\begin{equation} \label{e:definition_of_Psi_scheme}
\textstyle
\mathbf{\Psi} \colon \Gr_{\mathbb{Z}}(2, d)^r \dashrightarrow \Gr_{\mathbb{Z}}(r, \binom{d}{2}).
\end{equation}
We will tacitly assume throughout that $\mathbf{\Psi}$ is restricted to the open
subset on which it is regular.
Each Grassmannian scheme $\Gr_{\mathbb{Z}}(k, n)$ can be viewed as a 
projective scheme via Plücker coordinates. To be more precise,
consider the polynomial ring $P_k^n = \mathbb{Z}[..., X_I, ...]$
in $\binom{n}{k}$ variables indexed by the subsets $I \subseteq \{ 1, 2, \dots, n \}$
with $k$ elements. Consider the ring generated by matrix coefficients
$A_k^n = \mathbb{Z}[..., x_{i,j}, ...]$ of a generic matrix
$M = [x_{i,j}]_{1 \leq i \leq k, 1 \leq j \leq n}$.
To each $X_I$ we can associate the minor corresponding to columns indexed by $I$.
This gives a ring homomorphism $P_k^n \to A_k^n$ with kernel $J_k^n$. Then we have
\[
\Gr_{\mathbb{Z}}(k, n) = \Proj P_k^n / J_k^n
\subseteq 
\Proj P_k^n = \mathbb{P}^{\binom{n}{k} - 1}.
\]
Now, for the source scheme $\Gr_{\mathbb{Z}}(2,d)^r$, we will also need to use
the Segre embedding of the product of schemes into a projective space.
This corresponds to taking tensor products of rings, and we have
\[
\Gr_{\mathbb{Z}}(2,d)^r = \Proj \bigotimes_{i=1}^r P_2^d / J_2^d
\subseteq
\Proj \bigotimes_{i=1}^r P_2^d = \mathbb{P}^{\binom{d}{2}^r - 1}.
\]
Coordinates on the ring $\bigotimes_{i=1}^r P_2^d$ are indexed
as tuples $(I_1, \dots, I_r)$ with each $I_i$ a subset of $\{1,2, \dots, n\}$
of size $2$.
Fix a natural bijection between $\{ 1, 2, \dots, \binom{d}{2} \}$
and $\mathcal{I} = \{ I \subseteq \{ 1,2, \dots,d\} \mid |I| = 2 \}$.
Coordinates on the ring corresponding to $\Gr_{\ZZ}(r, W)$ can therefore
be indexed as $X_{\{ I_1, \dots, I_r \}}$ with $I_i \in \mathcal{I}$.
In this language, our map $\mathbf{\Psi}$ is given as a restriction
of the scheme version of the ring map
\[
P_r^{\binom{d}{2}} \to \bigotimes_{i=1}^r P_2^d,
\quad
X_{\{ I_1, I_2, \dots, I_r \}} \mapsto 
\sum_{\sigma \in \Sym_r} 
\mathrm{sgn}(\sigma) \cdot X_{( I_{\sigma(1)}, I_{\sigma(2)}, \dots, I_{\sigma(r)} )}
\]
Over fields, this corresponds to the projectivization of the natural projection
of $\otimes^r W$ onto $\wedge^r W$ symmetrizing a tensor.

\subsubsection*{Generic dimension of fibres}

Our main objective is to inspect the image of $\mathbf{\Psi}$ over
finite fields.
This will be done in two steps. The first step is to say something about
the generic behaviour of this map. We will do this by using
the following version of the fibre dimension theorem
(see \cite[\href{https://stacks.math.columbia.edu/tag/05F7}{Lemma 05F7}]{Stacks}).

\begin{lemma} \label{l:fiber_dimension_theorem}
Let $f \colon X \to Y$ be a morphism of schemes. 
Assume $Y$ irreducible with generic point $\eta$ and $f$ of finite type.
If $X_{\eta}$ has dimension $n$, then there exists a non-empty open $V \subseteq Y$
such that for all $y \in V$ the fibre $X_y$ has dimension $n$.
\end{lemma}

The dimension of the generic fibre $X_{\eta}$ will be obtained using submersiveness
based on the following.

\begin{lemma} \label{l:scheme_morphism_is_dominant}
Suppose $r \geq \binom{d-2}{2} + 1$. 
Then the rational map $\mathbf{\Psi}$ is dominant.
\end{lemma}
\proof
Since dominance is witnessed on open subsets (see 
\cite[\href{https://stacks.math.columbia.edu/tag/0CC1}{Lemma 0CC1}]{Stacks}),
we can replace $X$ and $Y$ by their
affine open subsets, and hence assume that $X = \Spec A$ and
$Y = \Spec B$ for finitely generated rings $A,B$ over
$\mathbb{Z}$.
Changing the base from $\mathbb{Z}$ to $\mathbb{C}$,
we have a diagram
\[
\begin{tikzcd}
X_{\mathbb{C}} \arrow[r, "\mathbf{\Psi}_{\mathbb{C}}"] \arrow[d, "\beta_X"] &
Y_{\mathbb{C}} \arrow[d, "\beta_Y"] \\
X_{\mathbb{Z}} \arrow[r, "\mathbf{\Psi}_{\mathbb{Z}}"] &
Y_{\mathbb{Z}},
\end{tikzcd}
\]
where $\beta_{\bullet}$ are base change projections induced by the scalar extensions
$A \to A \otimes_{\mathbb{Z}} \mathbb{C}$ and similarly for $B$.
Note that since $Y$ is irreducible and $Y_{\mathbb{Q}} \neq \emptyset$, 
the ring $B$ has no additive torsion. Therefore 
$B \to B \otimes_{\mathbb{Z}} \mathbb{C}$ is an injection, meaning
that $\beta_{Y}$ is dominant. It follows that the generic point of $Y_{\mathbb{Z}}$
belongs to the image of $\beta_Y$
(see \cite[\href{https://stacks.math.columbia.edu/tag/0CC1}{Lemma 0CC1}]{Stacks}). 
Now, by Theorem \ref{t:psi_submersion} the
image of $\mathbf{\Psi}_{\mathbb{C}}$ is not contained in any proper
subvariety, whence the rational map $\mathbf{\Psi}_{\mathbb{C}}$
is also dominant.
The generic point of $Y_{\mathbb{C}}$ is therefore in the image of
$\mathbf{\Psi}_{\mathbb{C}}$.
By commutativity of the diagram, we now have that the generic point of 
$Y_{\mathbb{Z}}$ is in the image of $\mathbf{\Psi}_{\mathbb{Z}}$.
This completes the proof.
\endproof

In order to express $\dim X_{\eta}$, we can pass to affine open subsets of both $X$
and $Y$ (see \cite[Exercise 3.20 (e)]{Hartshorne}) and hence assume that $X = \Spec A$ and $Y = \Spec B$ with a ring map
$B \to A$. This map is injective by Lemma \ref{l:scheme_morphism_is_dominant},
and so we can assume that $B \leq A$.
Therefore we have
\[
X_{\eta} =
\{ \mathfrak{p} \in \Spec A \mid \mathfrak{p} \cap B = 0 \} =
\Spec (B - \{ 0 \})^{-1} A.
\]
The latter corresponds to the spectrum of
$(B \otimes_{\mathbb{Z}} \mathbb{Q} - \{0\})^{-1} (A \otimes_{\mathbb{Z}} \mathbb{Q})$,
which is the same as the dimension of the generic fibre of 
$\mathbf{\Psi}_{\mathbb{Q}} \colon X_{\mathbb{Q}} \to Y_{\mathbb{Q}}$.
As this morphism is dominant and its source and target are irreducible, 
the latter dimension is equal to $\dim X_{\mathbb{Q}} - \dim Y_{\mathbb{Q}}$
(see \cite[Exercise 3.22 (b)]{Hartshorne}).
Whence we have
\[
\dim X_{\eta} = \dim X_{\mathbb{Q}} - \dim Y_{\mathbb{Q}}.
\]

%

\subsubsection*{Bounding rational points on fibres}

The second step in our reasoning will be to bound the sizes of the
fibres over finite fields. 
In order to achieve this, we will utilize the following bound, known
in the literature as the Schwarz-Zippel type bound (see \cite[Claim 7.2]{DKL14}).

\begin{lemma} \label{l:schwarz_zippel_bound}
Let $Z$ be a variety in $\overline{\mathbb{F}}_p^n$ defined over $\mathbb{F}_p$
of degree $d$ and dimension $k$. Then 
$| Z(\mathbb{F}_p) | \leq d \cdot p^k$.
\end{lemma}

We are now ready for the proof.

\proof[Proof of Theorem \ref{t:large_r_asymptotics}]

Set $X$ and $Y$ to denote the source and target of the rational map
$\mathbf{\Psi}$. 
We can replace $X$ and $Y$ by their affine open subsets,
and hence assume that both are affine varieties.
Let $\eta$ be the generic point of $Y$.
It follows from Lemma \ref{l:scheme_morphism_is_dominant}
that $\dim X_{\eta} = \dim X_{\mathbb{Q}} - \dim Y_{\mathbb{Q}}$.
Thanks to Lemma \ref{l:fiber_dimension_theorem}
there is a non-empty dense subset $V \subseteq Y$
containing $\eta$
such that $\dim \mathbf{\Psi}^{-1}(y) = \dim X_{\mathbb{Q}} - \dim Y_{\mathbb{Q}}$
for all $y \in V$. 
Set
$U = \mathbf{\Psi}^{-1}(V)$.
It follows from Theorem \ref{l:fiber_dimension_theorem}
applied to $U \to \Spec(\mathbb{Z})$
that for all but finitely many primes $p$,
we have $\dim U_{\mathbb{F}_p} = \dim U_{\mathbb{Q}}
= \dim X_{\mathbb{Q}}$.
Express
\[
\lvert U(\mathbb{F}_p) \rvert
=
\sum_{y \in \mathbf{\Psi}(U(\mathbb{F}_p))}
\lvert \mathbf{\Psi}^{-1}(y) \rvert.
\]
It follows from Lemma \ref{l:schwarz_zippel_bound} that
we can bound
\[
\lvert \mathbf{\Psi}^{-1}(y) \cap X(\mathbb{F}_p) \rvert
\leq
\deg \mathbf{\Psi}^{-1}(y)_{\overline{\mathbb{F}}_p} \cdot p^{\dim X_{\mathbb{Q}} - \dim Y_{\mathbb{Q}}}.
\]
Note that
$\deg \mathbf{\Psi}^{-1}(y)_{\overline{\mathbb{F}}_p} 
\leq
\deg X_{\overline{\mathbb{F}}_p}$
by B\'ezout (see~\cite[Chapter 8]{Ful98}) since $\mathbf{\Psi}$ is linear.
The latter degree is equal to
$\deg X_{\overline{\mathbb{Q}}} = (C_{d-2})^r$
(see \cite[Section 6]{Grassmannian}).
We thus obtain
\[
\lvert U(\mathbb{F}_p) \rvert
\leq
(C_{d-2})^r \cdot p^{\dim X_{\mathbb{Q}} - \dim Y_{\mathbb{Q}}}
\cdot
\lvert \mathbf{\Psi}(U(\mathbb{F}_p)) \rvert.
\]
The number of rational points of $X-U$
can be bounded by Lemma \ref{l:schwarz_zippel_bound}.
Together with Lemma \ref{l:grassmannian_precise_bound},
we obtain
\[
\lvert U(\mathbb{F}_p) \rvert 
\geq
p^{\dim X_{\mathbb{Q}}}
+ O(p^{\dim X_{\mathbb{Q}} - 1}),
\]
where the implied constant is independent of $p$.
It now follows that
\[
\lvert \im \psi \rvert
\geq
\lvert \mathbf{\Psi}(U(\mathbb{F}_p)) \rvert
\geq
\left( \frac{1}{C_{d-2}} \right)^r \cdot p^{\dim Y_{\mathbb{Q}}}
+
O(p^{\dim Y_{\mathbb{Q}} - 1})
\]
and the proof is complete.
\endproof

\begin{remark}

In general, it is not possible to replace the lower bound $1/(C_{d-2})^r$ with
$1$ for all values of $r$. This is clear from Proposition \ref{p:d=4}, where
the bound appears in the following way. A general point in $\Gr(2, W)$  is of
the form $\langle e_1 \wedge e_2, e_3 \wedge e_4 \rangle$ for a basis $\{ e_1,
e_2, e_3, e_4 \}$ of $V$. The fibre over such a point consists of the pair
$(e_1 \wedge e_2, e_3 \wedge e_4)$ and its permutation $(e_3 \wedge e_4, e_1
\wedge e_2)$. The generic fibre in this case consists of $(C_2)! = 2$ points,
while our method bounds these from above as $(C_2)^{r} = 4$. Proposition
\ref{p:d=4} indicates that the situation even improves with increasing $r$,
and not worsens as our bound does.

\end{remark}

\begin{remark}

Our method provides bounds that are interesting for increasing $p$ with $d,r$
fixed. Doing the opposite, i.e., fixing $p$ and increasing $d,r$, might result
in the large value of $(C_{d-2})^r$ as well as the implied constants
interacting with the numerical dimensions, and so the obtained bound will be
useless, possibly even negative for large enough $d,r$. We will show
in Theroem \ref{t:nontrivial_bog_1} that indeed in such a setting
it is possible to produce many elements outside $\im \psi$.

\end{remark}

\subsection{Log-generic $p$-groups} \label{ss:log_generic_groups}

Once the parameters $d$ and $r$ establishing our context are selected, each of the constructed groups $G_L$ will be of order
\[
|G_L| = |V| \cdot |W/L| = p^{d + \left( \binom{d}{2} - r \right)}.
\]
Now, rather than independently selecting
$d$ and $r$, suppose that we fix $n := d + \binom{d}{2} - r$ and with it the
sizes of the groups we are considering. 
Fix a parameter $0 < \alpha < 1$ and set $d := \lfloor \alpha n \rfloor$. 
Expressing $r = \binom{d}{2} + d - n$, we can therefore produce
groups of orders $p^n$ as $n$ varies and $\alpha$ is fixed.
The number of groups we obtain in this way is of the form
(see \cite{Hig60})
\begin{equation} \label{e:number_of_groups_alpha}
p^{\binom{d}{2} \left( \binom{d}{2} - r \right) + O(n^2)}
=
p^{\frac{\alpha^2}{2} (1 - \alpha) n^3 + O(n^2)},
\end{equation}
and since the size of $\GL_n(\FF_p)$ is at most $p^{O(n^2)}$,
the same expression gives an estimate on the number of
non-isomorphic groups we obtain.
In order to maximize the number of $p$-groups we obtain
in this way, one takes $\alpha = \frac23$, with which the above
simplifies into
\[
p^{\frac{2}{27} n^3 + O(n^2)}.
\]
It is well known (see \cite{Sim65})
that the latter also gives an asymptotical upper bound
for the number of {\em all} $p$-groups of order $p^n$.

We will now relate this to the previous section.
Using the above notation, we have $r - \binom{d-2}{2} = 3d - n - 4$,
which is asymptotically positive for values of $\alpha > \frac13$.
Thus for the optimal value $\alpha = \frac23$, 
we are in the situation of Theorem \ref{t:large_r_asymptotics}.
However, here we are fixing the prime $p$
and varying the parameters $d$ and $r$, whence the theorem does
not apply directly. We can instead use quite the opposite
Theorem \ref{p:small_r_asymptotics}. In order to maximize
its use, we should select $\alpha$ as large as possible
so that $\delta = \binom{d-2}{2} - r - 3$ is positive.
The latter is equivalent to $\alpha \leq \frac13$.
Thus we can parametrize groups with $\alpha = \frac13$,
and it follows from Theorem \ref{p:small_r_asymptotics}
that in this case, the number of elements
in the image of $\psi$ is at most $\frac12 |\Gr(r, W)|$.
Thus there are at least
$\frac12 |\Gr(r, W)| / |\GL(V)|$
non-isomorphic $p$-groups with non-trivial Bogomolov multipliers.
The size of the last Grassmannian can be bounded by
\[
\left| \Gr(r,W) \right| \geq p^{\left( \binom{d}{2} - 1 \right) r - r^2}
= p^{-\frac{d^3}{2}+\frac{d^2 n}{2} + O(n^2)}.
\]
As a function of $n$, the leading term of the exponential is $d^2(n-d)/2$,
which is of order $((\frac13 n)^2 \cdot \frac23 n) / 2 = \frac{1}{27}n^3$.
Factoring by the action of $\GL(V)$ does not change this order, since
the log-size of $\GL(V)$ is only quadratic in $n$. 
It therefore follows that as $n$ grows, we obtain many groups
with non-trivial Bogomolov multipliers, their number is of an order
of magnitude log-comparable to that of the number of all $p$-groups of order $p^n$.

The above is slightly unexpected, particularly so since we are in the
situation when 
$r \geq \binom{d-2}{2} + 1$
and so the corresponding $\QQ_p$-map is a local surjection
by Theorem \ref{t:psi_submersion},
meaning that many Bogomolov multipliers vanish.
The explanation is morally given 
by the fact that the constant in Theorem \ref{t:large_r_asymptotics}
converges to zero too quickly, so we can not obtain good
numerical bounds for the number of groups of order $p^n$ with vanishing
multipliers. 
At the same time, recall that by Proposition \ref{p:large_subvariety_Qp},
we can find a large subvariety avoiding the image of $\Psi$.
We now show how the same subvariety can be constructed over finite
fields at the value $\alpha = \frac23$.
This ultimately produces log-generic $p$-groups with 
{\em non}-vanishing Bogomolov multipliers.

Our construction is based on the following lemma. 
For a subspace $L \leq W$, set 
\[
L^{\wedge} := \langle D W \cap L \rangle =  \langle v \wedge w \mid v,w \in V, \; v \wedge w \in L \rangle.
\]

\begin{lemma} \label{l:not_in_image_of_psi}
Suppose $\rho := \binom{d}{2} - r \geq 3$. 
For every positive integer $N$, we have
\[
|\{ L \mid L \leq W, \; \dim L = r, \; \dim(L / L^{\wedge}) \geq N \}| \geq p^{(r - 3N)(\rho - 3N - 1)}
\]
as long as $d > 4N$, $r > 3N$ and $\rho > 3N + 1$. In particular,
taking $N=1$ gives
\[
|\Gr(r,W) - \im \psi|
\geq
p^{(r - 3) (\rho - 4)}.
\]
\end{lemma}
\proof

Let $\{ v_i \mid 1 \leq i \leq d \}$ be a basis of $V$ and let $\bar V_i \leq
V$ for $1 \leq i \leq N$ be the $4$-dimensional subspaces $\langle v_{4i - 3},
v_{4i - 2}, v_{4i - 1}, v_{4i} \rangle$. Set $X_i \leq \bar V_i \wedge \bar
V_i$ to be one of the subspaces of codimension $3$ that are not generated by
decomposable wedges (cf. the proof of Proposition \ref{p:d=4} (3)). 
Now let $U \leq W$ be the standard complement of
$\oplus_i ( \bar V_i \wedge \bar V_i )$, so that  $v_k \wedge v_l \in U$ as
long as not both $v_k, v_l$ belong to the same $\bar V_i$. Thus $W = \oplus_i
(\bar V_i \wedge \bar V_i) \oplus U$.

Select any subspace $L \leq U$ of dimension $r - 3N$. Now
\[
(\oplus_i X_i \oplus L)^{\wedge} \leq (\oplus_i X_i \oplus U)^{\wedge}
\]
and, as $U$ is generated by decomposable wedges, it follows that we have
a surjection 
\[
\frac
{\oplus_i X_i \oplus L}
{(\oplus_i X_i \oplus L)^{\wedge}}
\twoheadrightarrow
\frac
{\oplus_i X_i \oplus U}
{(\oplus_i X_i \oplus U)^{\wedge}}.
\]
The latter space is of dimension at least $N$, and therefore
\[
\dim (\oplus_i X_i \oplus L) / (\oplus_i X_i \oplus L)^{\wedge} \geq N.
\]
It follows that every subspace of $U$ of dimension $r-3N$ produces a 
subspace of $W$ of dimension $r$ whose decomposable wedges
belong to a subspace of codimension at least $N$. 
Thus the number of all such subspaces is at least
\[
\left| \Gr \left( r - 3N, U \right) \right| \geq
p^{(r - 3N)(\rho - 3N - 1)}. \qedhere
\]
\endproof

\begin{theorem} \label{t:nontrivial_bog_1}
Fix a prime $p > 2$ and a positive integer $M$.
Let $\#_{all}(n)$ be the number of all $p$-groups of order $p^n$, and let
$\#_{\B_0 \geq M}(n)$ be the number of $p$-groups of order $p^n$ whose Bogomolov
multiplier is of order at least $M$. Then
\[
\lim_{n \to \infty}
\frac
{\log_p \#_{\B_0 \geq M}(n)}
{\log_p \#_{all}(n)}
= 1.
\]
\end{theorem}
\proof
Set $\alpha = \frac23$ and consider the $p$-groups of order $p^n$ 
obtained with the above parametrization.
Thus
$d = \frac{2}{3} n + O(1)$,
$r = \frac{2}{9} n^2 + O(n)$.
Set $N = \log_p M$. Now, note that we have
\[
\#_{\B_0 \geq M}(n) \geq 
\frac
{|\{ L \mid L \leq W, \; \dim L = r, \; \dim(L / L^{\wedge}) \geq N \}|}
{|\GL(V)|}.
\]
Using Lemma \ref{l:not_in_image_of_psi} 
with $\rho = n - d = \frac{1}{3} n + O(1)$, it follows that
\[
\log_p \#_{\B_0 \geq M}(n) \geq r \rho + O(n^2) = 
\frac{2}{27} n^3 + O(n^2).
\]
The latter is also a log-upper bound for $\#_{all}(n)$
and the proof is complete.
\endproof

\subsection{Commutators in the end} \label{ss:commutators_in_the_end}

The methods in the previous section can be used in other types of problems
involving commutators. A sample application is the following elementary
statement about commutators in $p$-groups. Its proof captures the heart
of the argument of Theorem \ref{t:nontrivial_bog_1}.

\begin{proposition}
Fix a prime $p > 2$.
Let $\#_{all}(n)$ be the number of all $p$-groups of order $p^n$, and let 
$\#_{[G,G] \neq \K(G)}(n)$ be the number of $p$-groups of order $p^n$ 
in which not every element of the derived subgroup is a simple commutator. 
Then
\[
\lim_{n \to \infty}
\frac
{\log_p \#_{[G,G] \neq \K(G)}(n)}
{\log_p \#_{all}(n)}
= 1.
\]
\end{proposition}
\proof
As in the proof of Theorem \ref{t:nontrivial_bog_1},
consider only $p$-groups of order $p^n$ obtained by setting
$\alpha = \frac{2}{3}$.
Let $\bar V \leq V$ be a proper fixed subspace of dimension $4$
and let $U$ be the standard complement
of $\bar V \wedge \bar V$ in $W$ just like in the proof of Lemma \ref{l:not_in_image_of_psi}.
Every $r$-dimensional subspace $L \leq U$ gives a group $G_L$ of order $p^n$.
Factoring by the normal subgroup represented by $U$, we obtain a surjection
$G_L \to G_U$. The latter group is the free group of nilpotency class $2$
and exponent $p$ on $4$ generators, and not every element of its derived
subgroup is a simple commutator.
The latter property therefore also holds for the group $G_L$.
Thus we have
\[
\#_{[G,G] \neq \K(G)}(n) \geq \frac{|\Gr(r, U)|}{|\GL(V)|}.
\]
It follows that
\[
\log_p \#_{[G,G] \neq \K(G)}(n) \geq
r \left( \binom{d}{2} - r - 1 \right) + O(d^2) = 
\frac{2}{27} n^3 + O(n^2). \qedhere
\]
\endproof


\end{document}